\input amstex
\documentstyle{amsppt}
\input epsf
\NoBlackBoxes


\def\C{{\Bbb C}}
\def\R{{\Bbb R}}
\def\SS{{\Bbb S}}
\def\B{{\Bbb B}}
\def\P{{\Bbb P}}
\def\CP{{\Bbb CP}}
\def\Z{{\Bbb Z}}
\def\calU{{\Cal U}}
\def\calA{{\Cal A}}

\def\pr{\operatorname{pr}}
\def\supp{\operatorname{supp}}
\def\length{\operatorname{len}}

\def\Area{\operatorname{Area}}
\def\sign{\operatorname{sign}}
\def\dist{\operatorname{dist}}
\def\diam{\operatorname{diam}}
\def\const{\operatorname{const}}
\def\eps{\varepsilon}
\def\ph{\varphi}

\def\Re{\operatorname{Re}}
\def\Im{\operatorname{Im}}
\def\Arg{\operatorname{Arg}}

\def\pt{\operatorname{pt}}
\def\leg{\operatorname{leg}}
\def\card{\operatorname{Card}}
\def\<{\langle}
\def\>{\rangle}

\def\refBeloshapka{1}
\def\refProblems  {2}
\def\refBurglind  {3}

\def\sectLeg       {2}
\def\sectAppr      {3}
\def\sectStokes    {4}
\def\sectLegLB     {5}
\def\sectExplicit  {6}

\def\thMain       {1.1}
\def\remMiniMax   {1.2}
\def\remBidisk    {1.3}
\def\remMotiv     {1.4}
\def\remDegTwo    {1.5}

\def\lemLegOne    {2.1}
\def\lemLegTwo    {2.2}
\def\lemLegThree  {2.3}
\def\defNetwork   {2.4}
\def\defGenNet    {2.5}
\def\propLeg      {2.6}
\def\propGeneric  {2.7}

\def\lemApprOne   {3.1}
\def\remAppr      {3.2}
\def\corAppr      {3.3}
\def\defCircuit   {3.4}
\def\propAppr     {3.5}

\def\lemPerturbLin {3.6}
\def\lemPerturb    {3.7}
\def\propPerturb   {3.8}

\def\lemStokes    {4.1}
\def\lemDA        {4.2}
\def\defBetaPN    {4.3} 
\def\remBetaPN    {4.4}
\def\lemIsop      {4.5}
\def\remIsop      {4.6}
\def\lemLength    {4.7}
\def\corLength    {4.8}
\def\propEpsSq    {4.9}

\def\propLegLB    {5.1}

\def\propExplOne  {6.1}
\def\propExplTwo  {6.2}
\def\corExpl      {6.3}
\def\remExpl      {6.4}

\def\eqAppr        {1}
\def\eqApprOne     {2}
\def\eqTaylor      {3}
\def\eqApprTwo     {4}
\def\eqApprThree   {5}
\def\eqApprFour    {6}
\def\eqApprFive    {7}
\def\eqApprSix     {8}

\def\eqdrho        {10}

\def\eqNorm        {12}
\def\eqTwoPi       {13}
\def\eqStokes      {14}
\def\eqIsop        {15}
\def\eqDiam        {16}
\def\eqABL         {17}
\def\eqMain        {18}

\def\eqDiamLeg     {19}
\def\eqLegLBOne    {20}
\def\eqLegLBTwo    {21}
\def\eqLegLBThree  {22}
\def\eqLegLBFour   {23}

\def\eqEll         {24}
\def\eqSinCos      {25}
\def\eqExplOne     {26}
\def\eqExplTwo     {27}
\def\eqExplThree   {28}
\def\eqExplFour    {29}

\def\figLegNet        {1}
\def\figBetaPlus      {2}
\def\figBetaMinus     {3}
\def\figNet           {4}

\rightheadtext{ Algebraic curve in the unit ball in $\C^2$ }

\topmatter

\title    Algebraic curve in the unit ball in $\C^2$ passing
          through the center, all whose boundary components are
          arbitrarily short
\endtitle

\author   S.Yu.~Orevkov
\endauthor

\address  Universit\'e Paul Sabatier (Toulouse III)
\endaddress

\address  Steklov Math. Inst., Russian Acad. Sci.
\endaddress

\endtopmatter

\document

\centerline{\it To the memory of Anatoliy Georgievich Vitushkin }
\bigskip

\head  1. Introduction
\endhead

Let $\SS^3$ be the unit sphere in $\C^2$ centered in the origin.
A.G.~Vitushkin posed the following question (see
[\refBeloshapka], [\refProblems; Problem 5.3], [\refBurglind]):

\medskip
(1). 
Does there exist an absolute constant $c$ such that any complex algebraic
curve $A$ in $\C^2$ passing through the origin has a connected component
of the set  $A\cap\SS^3$ whose length is not greater than $c$?

(2). Is it true that $c=2\pi$?

\medskip

In this paper we give a negative answer to the both questions.

 \proclaim{ Theorem \thMain } 
 a).
Let $\Omega$ be a compact closed domain in an analytic surface and
$M=\partial\Omega$ its boundary. Let $M_0$ be the set of those points
where $M$ is a $C^2$-smooth strictly pseudoconvex real hypersurface.
Suppose that some Riemannian metric is fixed on $M$.
Let $A$ be a complex analytic curve in $\Omega$ such that $\partial A$ is contained
in $M_0$ and realizes the zero homology class in $H_1(M_0;\Z)$. Let $P$ be
any finite subset of $A$.

Then for any $2$-chain $\beta$ in $M_0$ such that $\partial\beta=\alpha$ and 
for any $\eps>0$, there exists a complex analytic curve $A'$ in $\Omega$
which is $\eps$-close to $A$ and such that the length of any its component
is less than $\eps$, and $P\subset A'$.

\smallskip

b). If, moreover, $\Omega\subset\C^2$ and, for any point $p\in M_0$, the complex
line $T$ tangent to $M$ at $p$ does not meet 
$\Omega$ at other points and the restriction to $T$ of the second fundamental form
of $M$ at $p$ is positive definite {\rm(}by the strict pseudoconvexity of $M$
the latter condition is equivalent to the fact that the sectional curvature at $p$
in the direction of $T$ is positive{\rm)} 
then one can choose $A'$ to be an algebraic curve.
\endproclaim

This theorem follows immediately from Propositions {\propLeg}
and {\propAppr}. It is proved in the end of \S\sectAppr.
The crucial role in the proof is played by the notion of a Legendrian net
hanged on a transversal cycle in a contact 3-manifold introduced in \S\sectLeg.

A negative answer to Vitushkin's question is provided by applying 
Theorem {\thMain} b) in the case when $\Omega$ is the unit ball, $P$ is its
center, and $A$ is an arbitrary curve (for example, a line) passing through $P$.

\remark{ Remark \remMiniMax }
We formulate Theorem {\thMain} and Propositions {\propLeg} and {\propAppr}
in a "minimax generality", i.e. we try to give a maximally general
statement under the condition that it can be proved using exactly the same
arguments as in the simplest (known to us) proof for the
case of an algebraic curve in the unit ball.

If one refuses of this principle then Theorem {\thMain} can be easily generalized
as far as one's fantasy allows. For instance, the line $T$ in Part b) could be replaced
by an algebraic curve (but then the proof of the corresponding analogue of Lemma
{\lemApprOne} would become more complicated), one could involve into consideration
Shilov boundaries, polynomial convexity, etc.
\endremark

\remark{ Remark \remBidisk } The condition of the {\sl strict} pseudoconvexity in
Theorem {\thMain} is important. Indeed, the answers to the both Vitushkin's questions
are positive if one considers the bidisk instead of the ball
(see [\refBeloshapka]).
\endremark

\remark{ Remark \remMotiv } Apparently, the main Vitushkin's motivation for asking
this question was its relation with the problems about polynomial hulls of "bad" sets.
Some links between these topics are discussed in a recent paper [\refBurglind].
\endremark

\remark{ Remark \remDegTwo } The answer to Question (2)
(is it true that $c=2\pi$) is negative even for curves of degree two.
To see this, one can explicitly parametrize the real curve
$\{\,(z,w)\in\C^2\,|\,w^2 = az(z-1)\,\}\cap\SS^3$, then check by numerical
integration that its length is less than $4\pi$ for some values of $a$, and
finally, to remark that the perturbed curve 
$\{w^2 = az(z-1+\eps)\}\cap\SS^3$ for $0<\eps\ll 1$ consists of two equal halves
whose total length is close to the length of the initial curve.
\endremark

\medskip

Thus, the absolute constant $c$ does not exist.
However if one fixes $n$ --- the number of connected components of
$A\cap\SS^3$, such a constant depending of $n$ certainly does exist
(it is clear that the total length of all the components is greater than $2\pi$).
Let $n(\eps)$ denote the minimal number of connected components of 
$A\cap\SS^3$ under the condition that $A$ is an algebraic curve through the
origin such that length of any connected component of 
$A\cap\SS^3$ is less than $\eps$.

It follows from the argument above that $n(\eps) > 2\pi/\eps$. It is not difficult
to deduce from Stokes' formula that after the projection onto $\CP^1$,
the sum of the oriented areas bounded by the projections of the components of
$A\cap\SS^3$ is greater than the area of the whole $\CP^1$, hence
$n(\eps) > \text{const}/\eps^2$ (see Proposition \propEpsSq).
On the other hand, a straight forward application of the construction
provided by the proof of Theorem {\thMain}, yields an upper bound
$n(\eps) < \text{const}/\eps^4$.

A natural correction of Vitushkin's question suggests itself:
is it true that the maximal length of components of $A\cap\SS^3$ 
is essentially greater than the evident estimates?
More precisely, what is the asymptotics of $n(\eps)$ as $\eps\to0$?
The same question can be asked about the quantity $d(\eps)$ --- the minimal
degree of an algebraic curve satisfying the same condition.
As we have seen, the order of growth of $n(\eps)$ is between $\eps^{-2}$
and $\eps^{-4}$. It seems plausible that it is $\eps^{-3}$.
In {\S\sectExplicit}, we prove an upper bound for
$n(\eps)$ of the order $\eps^{-3}$.
In {\S\sectLegLB}, we prove that this bound cannot be improved by the methods
of this paper (i.e. using the construction based on a perturbation of a
Legendrian net). 
In the end of {\S\sectLegLB}, we propose a new question, a positive answer to whom
would imply a lower bound on  $n(\eps)$ of the order $\eps^{-3}$.

%
%

\head \sectLeg. Legendrian nets hanged on transverse cycles
\endhead

All the statements of this section are almost evident but we shall give
however their proofs.

We shall understand
{\it chains, cycles}, and {\it boundaries} more or less in the sense of
the theory of singular homologies, but we shall consider only piecewise
smooth chains and we shall identify chains obtained one from another by 
subdivisions and reparametrizations. In particular,
a 1-{\it chain} in a smooth manifold $M$ is by definition an element
of the quotient
of the free abelian group generated by all piecewise smooth mappings
$\alpha:I=[0,1]\to M$ modulo all relations of the form
$\alpha=-(\alpha\circ\varphi)$ and
$\alpha=(\alpha\circ\varphi_1) + (\alpha\circ\varphi_2)$ where $\varphi$ is
an orientation reversing piecewise smooth homeomorphism of the segment $I$
onto itself, and $\varphi_1,\varphi_2$ are orientation preserving 
piecewise smooth homeomorphisms of the segment $I$
onto the segments  $[0,1/2]$ and $[1/2,1]$ respectively.
For example, these relations imply that the constant mapping 
$I\to p\in M$ realizes the zero chain.
A linear combination $\sum m_i\alpha_i$ representing a chain $\alpha$
will be called a {\it minimal realization of} $\alpha$ if 
$m_i\ne 0$ for all $i$ and there does not exist indices $i_1$, $i_2$,
segments $I_1,I_2\subset I$, and a homeomorphism $\varphi:I_1\to I_2$ such that
$\alpha_{i_1}|_{I_1} = \alpha_{i_2}\circ\varphi$.

Let $\sum m_i\alpha_i$ be some minimal realization of a $1$-chain 
$\alpha$ on a 3-manifold $M$. Then the set $\supp\alpha=\bigcup\alpha_i(I)$
is called the {\it support} of $\alpha$. If, moreover, $M$ is endowed with a
Riemannian metric then the {\it length} of $\alpha$ is by definition 
$\length\alpha = \sum |m_i|\length\alpha_i$ where 
$\length\alpha_i$ is the length of the path $\alpha_i$.
A 1-chain  $\alpha$ is called
{\it $\eps$-short} if $\length\alpha<\eps$.
Analogously we define the support and the area of a 2-chain.
In the sequel, we shall not distinguish between chains and their minimal realizations.
A 1-cycle is called {\it generic} or {\it in general position} if it is a union
of pairwise disjoint piecewise smoothly embedded oriented circles taken with the
multiplicity 1.

Recall that a {\it contact structure} on a 3-manifold $M$ is a smooth field of 2-planes
which can be represented as $\ker\eta$ where $\eta$ is a 1-form such that
$\eta\wedge d\eta$ does not vanish.
It is known that all contact structures are locally equivalent to each other.

A 1-chain on a contact 3-manifold $(M,\eta)$ is called
{\it Legendrian} if it is $C^2$-smooth and the restriction of $\eta$ identically
equals to zero on its smooth pieces.
A 1-chain $\alpha$ is called {\it positively transverse} if it can be represented
as $\alpha=\sum m_i\alpha_i$ where $m_i>0$ and
$\alpha_i^*(\eta)>0$ for all $i$ (such a realization of $\alpha$
automatically is minimal).


Let us denote the standard coordinates in $\R^3$ by $x,y,z$ and let us 
consider the contact structure defined by the 1-form $\eta = dz-y\,dx$.
Let $\pr:\R^3\to\R^2$ be the projection $(x,y,z)\mapsto(x,y)$.

\proclaim{ Lemma \lemLegOne } Let $\gamma:[0,1]\to\R^2$ be a $C^2$-smooth path
starting at a point $p_0=(x_0,y_0)$. Then for any 
$z_0\in\R$ there exists a unique Legendrian path
$\tilde\gamma:[0,1]\to\R^3$ starting at $\tilde p_0=(x_0,y_0,z_0)$ such that
$\gamma = \pr\tilde\gamma$. Moreover, the length of
$\tilde\gamma$ is less than $L\sqrt{1+(|y_0|+L)^2}$ where $L$ is the length of
$\gamma$.
\endproclaim

The path $\tilde\gamma$ is called the {\it Legendrian lift} of
$\gamma$ starting at $\tilde p_0$.

\demo{ Proof } Let $\gamma(t)=(x(t),y(t))$. Set
$\tilde\gamma(t) = (x(t),y(t),z(t))$ where
$z(t) = z_0 + \int_{\gamma([0,t])} y\,dx$.
We have $|\tilde\gamma'|^2 = \dot x^2 + \dot y^2 + \dot z^2
=\dot x^2 + \dot y^2 + (y\,\dot x)^2 \le (1+y^2)|\gamma'|^2$.
Hence, the length of $\tilde\gamma$ is less than
$L\max\sqrt{1+y^2}$. It remains to note that $\max y \le |y_0|+L$.
\qed\enddemo

\proclaim{ Lemma \lemLegTwo } Let $0<\eps<1/2$ and let $\tilde
p_0=(x_0,y_0,z_0)$ and $\tilde p_1=(x_1,y_1,z_1)$ be points in
$\R^3$ such that $|y_0|<1$, $|y_1|<1$, and $\|\tilde p_1-\tilde
p_0\|<\eps^2$. Then there exists a piecewise smooth Legendrian path from
$\tilde p_0$ to $\tilde p_1$ whose length is less than $c_1\eps$
for some absolute constant $c_1$.
\endproclaim

\demo{ Proof } Let $\gamma_1$ be a straight line segment connecting
$p_0=\pr(\tilde p_0)$ to $p_1=\pr(\tilde p_1)$ and let
$\tilde\gamma_1$ be the Legendrian lift of $\gamma_1$ starting at
$\tilde p_0$. Let $\tilde p_1'=(x_1,y_1,z'_1)$ be the end of
$\tilde\gamma_1$. Let $\gamma_2=\sign(z_1-z'_1)\partial D$ where
$D$ is a disk of area $|z_1 - z'_1|$ such that $p_1\in\partial D$.
Let $\tilde\gamma_2$ be the Legendrian lift of $\tilde\gamma_2$ starting at
$\tilde p_1'$. Then the end of $\tilde\gamma_2$ coincides with
$p_1$ because $\int_{\tilde\gamma_2}dz = \int_{\gamma_2}y\,dx =
\pm\Area(D)$. The estimate for the length of $\tilde\gamma_1$ 
is obtained by the straight forward application of Lemma \lemLegOne. \qed
\enddemo

\proclaim{ Lemma \lemLegThree } Let $M$ be a contact $C^2$-smooth $3$-manifold
endowed with a Riemannian metric. Let
$\alpha$ be a Legendrian zero-homologous $1$-cycle on $M$.
Then for any $\eps>0$ there exist $\eps$-short Legendrian $1$-cycles
 $\alpha_1,\dots,\alpha_n$ on $M$ such that $\sum\alpha_j=\alpha$.
\endproclaim

\demo{ Proof } 
It is known that all contact structures are locally equivalent to each other.
Hence, for any $p\in M$ there exist its neighbourhood
$U_p$ and a smooth embedding $\ph_p:U_p\to\R^3$ taking the given contact 
structure on  $M$ to the contact structure on $\R^3$ defined by
the form
$\eta=dz-y\,dx$. Replacing if necessarily $U_p$ by a smaller neighbourhood,
we may assume that the set $\ph_p(U_p)$ is convex, contained in the layer
$\{|z|<1\}$, and there exists a constant
$m_p>0$ such that $\|d\ph_p(v)\|>m_p\|v\|$ for all $v\in TU_p$.
In each $U_p$, let us choose an open subset $V_p$ such that $p\in
V_p$ and $\overline V_p\subset U_p$.

Let $\beta$ be a $2$-cycle in $M$ whose boundary is $\alpha$.
Let us choose a finite subfamily $\calU=\{(U_i,V_i,\ph_i)\}_{i=1,\dots,k}
\subset\{(U_p,V_p,\ph_p)\}_{p\in M}$
such that the support of $\beta$ is contained in $\bigcup_{i=1}^k V_i$\,,
and let $m=\min_{(U_p,V_p,\ph_p)\in\calU}m_p$.

Let $\eps_1 = \min_i \dist\big( \ph_i(\overline V_i),\,\R^3\setminus\ph_i(U_i)\big)$
and let $\eps_2 = \min(\eps_1, m\eps/3)/c_1$
(here $c_1$ is the constant from Lemma \lemLegTwo).
Let us represent $\beta$ as a sum of simplices
$\beta=\beta_1+\dots+\beta_n$ so that:
\roster
\item
   each $\beta_j$ is contained in some $V_{i_j}$
   and $\diam\ph_{i_j}(\beta_j) < \eps_2^2$ far any $i=1,\dots,n$;
\item
   the lengths (with respect to the metric on $M$)
   of those edges of $\beta_j$'s which contribute to
   $\alpha=\partial\beta$ are smaller than $\eps/3$.
\endroster

Let $\Gamma=\{\gamma\}$ be the set of those edges of the simplices
$\beta_j$'s which do not contribute to $\alpha$ (for each pair of
edges which cancel against each other in $\partial\sum\beta_j$, we
include only one of them to $\Gamma$). For each $\gamma\in\Gamma$,
$\gamma\subset\partial\beta_j$, using Lemma \lemLegTwo, we can
choose a piecewise Legendrian path $\gamma'$ which relates the
ends of $\ph_{i_j}(\gamma)$ and which is shorter than $c_1\eps_2$.
Since $c_1\eps_2\le\eps_1$, we have $\gamma'\subset U_{i_j}$,
hence $\gamma''=\ph_{i_j}^{-1}(\gamma')$ is a Legendrian path in
$M$ shorter than $\eps/3$. Let $\Gamma''$ be the set of all those
$\gamma''$.

Finally, for each $j=1,\dots,n$, we define $\alpha_j$ as the cycle
obtained from the boundary of $\beta_j$ by replacing every its
edge $\gamma\in\Gamma$ with the corresponding path
$\gamma''\in\Gamma''$. \qed
\enddemo

\definition{ Definition \defNetwork }
Let $\alpha$ be a positively transverse 1-cycle in a contact
3-manifold $M$. A finite collection of 1-cycles
$\calA=\{\alpha_1,\dots,\alpha_n\}$ in $M$ is called a 
{\it Legendrian net hanged on} $\alpha$ if
\roster 
\item
       each $\alpha_i$ decomposes into the sum of two 1-chains
       $\alpha_i=\alpha_i^{\pt}+\alpha_i^{\leg}$
       (each of them may be zero)
       where $\alpha_i^{\pt}$ is positively transverse and
       $\alpha_i^{\leg}$ is Legendrian;
\item
       $\alpha_1+\dots+\alpha_n=\alpha^{\pt}_1+\dots+\alpha^{\pt}_n=\alpha$;
\endroster
The cycles $\alpha_1,\dots,\alpha_n$ will be called the {\it cells} of $\calA$,
and the union of their supports will be called the {\it support} of $\calA$.
\enddefinition

\definition{ Definition \defGenNet }
 Let $\alpha$ be a generic positively transverse cycle in $M$.
 A Legendrian net $\calA=\{\alpha_1,\dots,\alpha_n\}$ hanged on
 $\alpha$ is called {\it generic} or {\it in generic position}
 if there exists a piecewise smoothly embedded graph $\Gamma$
 with Legendrian edges such that
\roster
\item
      the multiplicity (i.e. the number of incident edges) of any vertex of $\Gamma$
      is either 1 or 3;
\item
      each end (i.e. a vertex of multiplicity 1) of $\Gamma$ is a smooth point 
      of the support of $\alpha$, and the tangents to $\Gamma$ and to $\alpha$ at
      this point are distinct;
\item
      $\Gamma\cap\supp\alpha$ coincides with the set of the ends of $\Gamma$;
\item
      each chain $\alpha_i^{\leg}$ is a sum of edges of $\Gamma$ taken with 
      the coefficients $\pm1$, every edge contributing to exactly two cells with
      the opposite signs.
\endroster

\enddefinition

\proclaim{ Proposition \propLeg } Let $M$ be a contact
$C^2$-smooth $3$-manifold endowed with a Riemannian metric.
Let $\alpha$ be a positively transverse $1$-cycle in $M$ which is
homologous to zero. Then for any $\eps$ there exists a Legendrian net
$\calA=\{\alpha_1,\dots,\alpha_n\}$ в hanged on
$\alpha$ all whose cells are $\eps$-short.
The support of $\calA$ can be done arbitrarily close to the support of
any given $2$-chain $\beta$ such that
$\partial\beta=\alpha$.

If, moreover, $\alpha$ is generic then $\calA$ also can be done generic.
\endproclaim

The proof is the same as that of  {\lemLegThree}, and we omit it.
To achieve the genericity of $\calA$, one should apply the following statement.

\proclaim{ Proposition \propGeneric }
Let $M$ be a contact
$C^2$-smooth $3$-manifold endowed with a Riemannian metric.
Let $\alpha$ be a generic positively transverse $1$-cycle in $M$, and let
$\calA=\{\alpha_1,\dots,\alpha_n\}$ be a Legendrian net hanged on $\alpha$.

Then for any $\delta>0$ there exists a generic Legendrian net 
$\calA'=\{\alpha'_1,\dots,\alpha'_n\}$ hanged on $\alpha$ such that
for any $i=1,\dots,n$, the cycles $\alpha'_i$ and $\alpha_i$ are $\delta$-close
in Hausdorff metric and $|\length\alpha_i - \length\alpha'_i|<\delta$.
\endproclaim

\demo{ Proof }
{\it Step 1.} Let us show that after an arbitrarily small perturbation of $\calA$,
one can find an embedded graph $\Gamma$ with Legendrian edges such that 
Condition (4) of Definition {\defGenNet} is satisfied.

\smallskip
By the definition of $1$-chains, there exist piecewise smooth Legendrian paths
$\gamma_1,\dots,\gamma_k$ and integer coefficients $m_{ij}$ such that
$\alpha_i^{\leg}=\sum_j m_{ij}\gamma_j$.
We must achieve $|m_{ij}|\le1$ for all $i,j$. 
To this end we shall successively reduce the quantity 
\newline $\sum_{ij}\max(0,|m_{ij}|-1)$.
Suppose that $m_{i_0,j_0}\ge2$ for some $i_0$, $j_0$
(the case $m_{{i_0},{j_0}}\le-2$ is analogous).
Since $\gamma_{j_0}$ does not contribute to $\sum_i\alpha_i$, we have
$\sum_i m_{i,j_0}=0$. Hence, there exists an index $i_1$ such that
$m_{i_1,j_0}<0$. Let $\gamma'$ be a Legendrian perturbation of 
$\gamma_{j_0}$ such that $\partial\gamma'=\partial\gamma_{j_0}$ and
$(\supp\gamma')\cap(\supp\Gamma)=\supp(\partial\gamma')$.
Let us replace $\alpha_{i_0}$ with $\alpha_{i_0}-\gamma_{j_0}+\gamma'$
and $\alpha_{i_1}$ with $\alpha_{i_1}+\gamma_{j_0}-\gamma'$.
It easy to see that this reduces the quantity $\sum_{ij}\max(0,|m_{ij}|-1)$
at least by one.

\medskip
{\it Step 2.} Suppose that there exists an embedded graph $\Gamma$ with 
Legendrian edges which satisfies Condition (4) of Definition \defGenNet, 
and let us show that it can be perturbed so that (1)--(3) are satisfied.
\smallskip

Let $p$ be a vertex of $\hat\Gamma=\Gamma\cup(\supp\alpha)$
of multiplicity $k>3$. Let us consider an auxiliary graph $G_p$ defined as follows.
Its vertices are the edges of $\hat\Gamma$ incident to $p$. Two
vertices $\gamma$ and $\gamma'$ of $G_p$ (i.e. edges of $\hat\Gamma$) are
connected by an edge in $G_p$ when $\gamma\subset\supp\alpha_i$ and
$\gamma'\subset\supp\alpha_i$ for some $\alpha_i$.
The condition $\sum_i\alpha_i=\alpha$ implies that after removing certain edges
from $G_p$, one obtains a disjoint union of graphs
$E_1\sqcup \dots\sqcup E_m$, moreover, the graphs $E_2,\dots,E_m$
are combinatorially equivalent to a circle, and $E_1$ is equivalent either to a circle
(when $p\not\in\supp\alpha$), or to a segment whose endpoints correspond
to the edges of $\hat\Gamma$ lying on $\supp\alpha$.
Denote the vertices of $E_k$ by
$\gamma_{k,1},\dots,\gamma_{k,c_k}$ so that $\gamma_{k,j}$ is connected to 
$\gamma_{k,j+1}$ by an edge in $E_k$.

Let $U_p$ be a sufficiently small neighbourhood of $p$ diffeomorphic to the ball,
such that 
$\Gamma_p=U_p\cap\Gamma=\bigcup_{k,j}(\gamma_{k,j}\cap U_p)$ and each of
$\gamma_{k,j}\cap U_p$ is an embedded segment transverse to $\partial U_p$. Set 
$\Gamma_{p,k} = \bigcup_{j=1}^{c_k}(\gamma_{k,j}\cap U_p)$ and 
$q_{k,j}=\gamma_{k,j}\cap\partial U_p$.
Let $\Gamma'_{p,k}$, $k=1,\dots,m$, be an arbitrary plane tree embedded into
a disk $\Delta$ all whose vertices having the multiplicity 1 or 3,
the number of ends (i.e. vertices of multiplicity 1) being equal to $c_k$,
and all the ends lying on $\partial\Delta$.
Let us denote the ends of $\Gamma'_{k,j}$ by
$q'_{k,1},\dots,q'_{k,c_k}$ in this cyclic order along $\partial\Delta$.
When $p\in\supp\alpha$, we shall also assume that there exists a vertex $p'$ of
$\Gamma'_{p,1}$ connected by edges to $q'_{1,1}$
and $q'_{1,c_1}$.

To perturb $\Gamma$ as is required, we replace each tree $\Gamma_{p,k}$ by the image of
$\Gamma'_{p,k}$ under an embedding into $M$ which has the following properties.
It takes $q'_{k,j}$ to $q_{k,j}$, it maps homeomorphically the union of the edges
$[p',q'_{1,1}]\cup [p',q'_{1,c_1}]$ onto the arc 
$q_{1,1}q_{1,c_1}$ of $\alpha$
(the vertex $p'$ being sent to a smooth point of this arc), and the
images of all other edges of $\Gamma'_{p,k}$ are Legendrian.
\qed
\enddemo

%
%

\head \sectAppr. Approximation of a Legendrian net by the union of 
                 boundaries of analytic disks
\endhead

Let $V$ be a complex analytic surface, and $M$ a real hypersurface in $V$.
Then the field of complex tangents is defined on $M$. It can be represented as
$\ker\eta$ for some $1$-form $\eta$. We shall call a curve $\gamma:[0,1]\to M$
{\it Legendrian} (resp. {\it positively transverse}) if $\gamma^*\eta=0$
(resp. $\gamma^*\eta>0$). In the case when $M$ is strictly pseudoconvex,
the field of complex tangents is a contact structure on $M$, hence these 
definitions are coherent with the definitions in \S\sectLeg.

\proclaim{ Lemma \lemApprOne } Let $U$ be an open subset in 
$\C^2$, and $M\subset U$ a real hypersurface defined by an equation
$f=0$ where $f$ is a real $C^2$-smooth function in $U$.
Let $\gamma:[0,t_1]\to M$ be a Legendrian $C^2$-smooth path and let
$p_0=\gamma(0)$. Let $T$ be the complex tangent line to 
$M$ at $p$. Suppose that the Hessian $H$ at $p_0$ of the restriction
$f|_T$ is positive definite. 

Let $L_t$ be the complex line passing through the points $p$ and $\gamma(t)$.
Let $S_t^+$ and $S_t^-$ denote the arcs into which the curve $L_t\cap M$
is divided by the points $\gamma(0)$ and $\gamma(t)$.
Then we have
$$
       \lim_{t\to 0}\; {2\length(S_t^\pm)\over\length\big(\gamma([0,t])\big)}
       = {\length(E)\over d(E,\gamma'(0))}
       < \pi\sqrt{K_1/K_2},                                        \eqno(\eqAppr)
$$
where $E$ is the ellipse $\{H=1\}$, $d(E,v)$ is the length of its diameter in the
direction of a vector $v$, and 
$K_1$, $K_2$  $\;(K_1\ge K_2)$ are the principal curvatures of $M$
in the direction of $T$.
\endproclaim

\demo{ Proof }  Let us denote the coordinates in $\C^2$ by $(z,w)$.
Without loss of generality we may assume that $p_0$ is the origin,
$T$ is the axis $w=0$, and $\gamma'(0)=(1,0)$. Then we have
$$
     f'_z(0,0)=f'_{\bar z}(0,0)=0
     \qquad\text{and}\qquad
     f'_w(0,0) = \overline{f'_{\bar w}(0,0)} = a \ne 0.    \eqno(\eqApprOne)
$$
Since $f$ is twice differentiable, we have
$$
     f(z,w) = aw + \bar a\bar w + Az^2 + 2Bz\bar z + \bar A\bar z^2
                 + w\,g_1 + \bar w\,g_2 + (z\bar z+w\bar w)\,g_3,
                                                           \eqno(\eqTaylor)
$$
where
$$
    2A=f''_{zz}(0,0),\quad
    2B=f''_{z\bar z}(0,0),\quad
    \lim_{(z,w)\to(0,0)}g_{1,2,3}(z,w)=0.
$$
Let us set $\gamma(t) = (z(t),w(t))$. The condition that the path $\gamma$
is Legendrian means that 
$$
 f'_z(\gamma(t))\,z'(t) + f'_w(\gamma(t))\,w'(t) = 0,\qquad
        t\in[0,t_1].                                      \eqno(\eqApprTwo)
$$
For $t=0$, by (\eqApprOne), this implies $w'(0)=0$. Hence we have
$$
        z(t) = t(1+\alpha_1(t)),\quad
        w(t) = b\,t^2(1+\alpha_2(t)),\quad
        2b=w''(0),\quad
        \lim_{t\to 0}\alpha_{1,2}(t) = 0.                \eqno(\eqApprThree)
$$
Differentiating  (\eqApprTwo) at $t=0$ and combining with
(\eqApprOne), (\eqTaylor), and (\eqApprThree), we obtain
$$
      2 ab + 2A + 2B = 0.
                                                         \eqno(\eqApprFour)
$$
Consider the parametrization of $L_t$ given by
$\ph_t:\C\to\C^2$, $\zeta\mapsto\big(z(t)\zeta,w(t)\zeta\big)$.
Let us denote the curve $\ph_t^{-1}(L_t\cap M)$ by $S_t$. It is defined by
$f(\ph_t(\zeta))=0$. Using
(\eqTaylor) and (\eqApprThree), one can rewrite the left hand side of this 
equation in the form
$$
    t^2\cdot\Big( ab\zeta + \bar a\bar b\bar\zeta + A\zeta^2
    + 2B\zeta\bar\zeta + \bar A\bar\zeta^2 + g(t,\zeta)\Big)
                                                         \eqno(\eqApprFive)
$$
where $g(t,\zeta)$ tends to zero as $t\to0$ uniformly on any bounded subset of
$\C$. Note, that the Hessian of the restriction 
$f|_T$ has the form $H(z)=(Az^2 + 2Bz\bar z + \bar A\bar z^2)/2$. Hence,
combining (\eqApprFive) with (\eqApprFour) and dividing by $2t^2$,
we obtain $S_t = \{\,\zeta\,|\,
H(\zeta-1/2)+g(t,\zeta)=H(1/2)\,\}$, and hence, $S_t\to E_{1/2}$ for
$t\to0$ where $E_{1/2} = \{\,\zeta\,|\,H(\zeta-1/2)=H(1/2)\,\}$
is a translate of $E$. Since the second derivatives of $f$
are continuous,
$S_t\to E_{1/2}$ implies $\length(S_t)\to\length(E)$ and
$\length(S_t^\pm)\to\length(E)/2$. 
It remains to note that the length of
$\ph_t^{-1}(\gamma[0,t])$ tends to 
$d(E,1)$, because $\gamma$ is twice differentiable.
\qed\enddemo

\remark{ Remark \remAppr } a). If, in Lemma {\lemApprOne}
 we replace the condition that $\gamma$ is
Legendrian by a weaker condition that $\gamma'(0)\in T$, then
$S_t$ would still tend to some translate of $E$. However, it may happen
in this case, that the center of the translated ellipse would not be on 
the real axis, hence the arc $\ph_t^{-1}(\gamma[0,t])$ would not tend to
a diameter. Thus, the upper bound for the ratio of the length would fail.

\smallskip
b). The only place in the proof where the continuity of the second derivatives
of $f$ is used, is the implication
$(S_t\to E_{1/2}) \Longrightarrow(\length(S_t)\to\length(E))$.
Therefore, the assertion of the lemma remains true if we replace the condition
that $f$ is of the class $C^2$ by a weaker condition that $f$ is just twice
differentiable, but if we assume in addition that $M$ is convex.
\endremark

\proclaim{ Corollary \corAppr }
Let $\Omega$ be a domain in a complex analytic surface whose
boundary $M=\partial\Omega$ is $C^2$-smooth. Suppose that $M$ is endowed with
a $C^2$-smooth Riemannian metric $g$, and let $p_0\in M$.
Suppose that $M$ is strictly pseudoconvex in a neighbourhood of a point $p_0$.
Let $\gamma:[0,t_1]\to M$, $\gamma(0)=p_0$,
be a Legendrian $C^2$-smooth curve.

Then, for any $\delta>0$, there exists a family of analytic disks
$\{D_t\}_{t\in[0,t_2]}$, $t_2\le t_1$ such that
$D_t\subset\Omega$, $\partial D_t\subset M$,
$D_t\cap\gamma = \{p_0,\gamma(t)\}$, $D_t$ is transverse to $M$,
and
$$
       \lim_{t\to 0}\; {\length S_t^\pm\over\length\gamma([0,t])}
       < 1+\delta
                                                        \eqno(\eqApprSix)
$$
where $S_t^+$ and $S_t^-$ are the arcs into which the curve
$\partial D_t$ is divided by the points $\gamma(0)$ and $\gamma(t)$.
\endproclaim

\demo{ Proof }
Let us choose the coordinates $(z,w)$ as in the proof of Lemma {\lemApprOne}.
Then the coordinate change $(z,w)\to(z,w+cz)$
transforms $H$ into
$$
     (A+ac)z^2 + 2Bz\bar z + (\bar A+\bar a\bar c)\bar z^2.
$$
Let us choose $c$ so that $A+ac=B-\delta_1$ when $\delta_1\ll\delta$,
and apply Lemma {\lemApprOne}.
\qed\enddemo

\definition{ Definition \defCircuit }
Let $M$ be a smooth contact manifold. A {\it Positive Transverse Simple
Crossing Curve} (PTSC-curve) 
on $M$ is a union of piecewise smooth embedded positively transverse oriented 
closed curves $S=S_1\cup\dots\cup S_m$ (called the {\it components of} $S$)
which meet each other at most pairwise and so that
if $S_i$ and $S_j$ intersect at $p$ then each of these curves is smooth at $p$
and the tangents to $S_i$ and to $S_j$ at $p$ are distinct.

A {\it circuit} of a PTSC-curve $S$ is an oriented piecewise smooth embedded circle
$\gamma$ which is a union of arcs of
$S$ such that 
\roster
\item
      on any smooth arc $a$ of $\gamma$, the orientation induced from $\gamma$
      coincides with the orientation induced from $S$;
\item
      if $\gamma$ passes through the intersection point of two
      components of $S$ then it switches from one component to the other one.
 \endroster

It is clear that any two circuits may intersect each other only at 
intersection points of components of $S$, and the sum of all circuits is $S$.
\enddefinition

\proclaim{ Proposition \propAppr }
 a). 
Let $\Omega$ be a domain in a complex analytic surface whose
boundary $M=\partial\Omega$ is $C^2$-smooth. Suppose that $M$ is endowed with
a $C^2$-smooth Riemannian metric $g$.
Let 
$\alpha$ be a positively transverse curve which is a union of disjoint piecewise
$C^2$-smoothly embedded circles, and let
$\calA=\{\alpha_1,\dots,\alpha_n\}$ be a generic Legendrian net hanged on 
$\alpha$. Suppose that $M$ is strictly pseudoconvex at a neighbourhood of 
$\supp\calA$.

Then, for any $\delta>0$, there exists a PTSC-curve $S=S_1+\dots+S_N$
such that:
\roster
\item
     each $S_j$ is the boundary of an analytic disk $D_j$ in $\Omega$;
\item
     $S+\alpha$ has exactly $n$ circuits $\beta_1,\dots,\beta_n$;
\item
     for any $i=1,\dots,n$, the Hausdorff distance between 
     $\beta_i$ and $\alpha_i$ is less than $\delta$, and
     $\length(\beta_i) < \length(\alpha_i)+\delta$.
\endroster

\medskip
b). If, moreover, $\Omega$ is a domain in 
$\C^2$ and the sectional curvature of $M$ in the direction of complex
tangents does not vanish in some neighbourhood of $M$, then
the disks $D_1,\dots,D_N$ can be chosen so that each of them is the intersection
of $\Omega$ with some complex line, but in this case the estimate for the lengths
should be replaced by 
$\length(\beta_i) < c_3\length(\alpha_i) + \delta$
where $c_3$ is a constant depending on $M$, $g$, and $\calA$.
{\rm(}when $\Omega$ is the unit ball and $g$ is induced by the standard metric in
$\C^2$, one has $c_3=\pi/2${\rm)}.
\endproclaim

\demo{ Proof } 
a).
Induction by $n$. The case $n=0$ is trivial.
Suppose that we proved the required statement for Legendrian nets having
$n-1$ cells. Let us prove it for a Legendrian net $\calA$ which has $n$ cells. 
Let $\alpha_i^{\pt}$ and $\alpha_i^{\leg}$ be as in Definition \defNetwork.
The construction described below is illustrated in Figures 
\figLegNet(a--d).

\midinsert
\centerline{
\epsfxsize 55mm\epsfbox{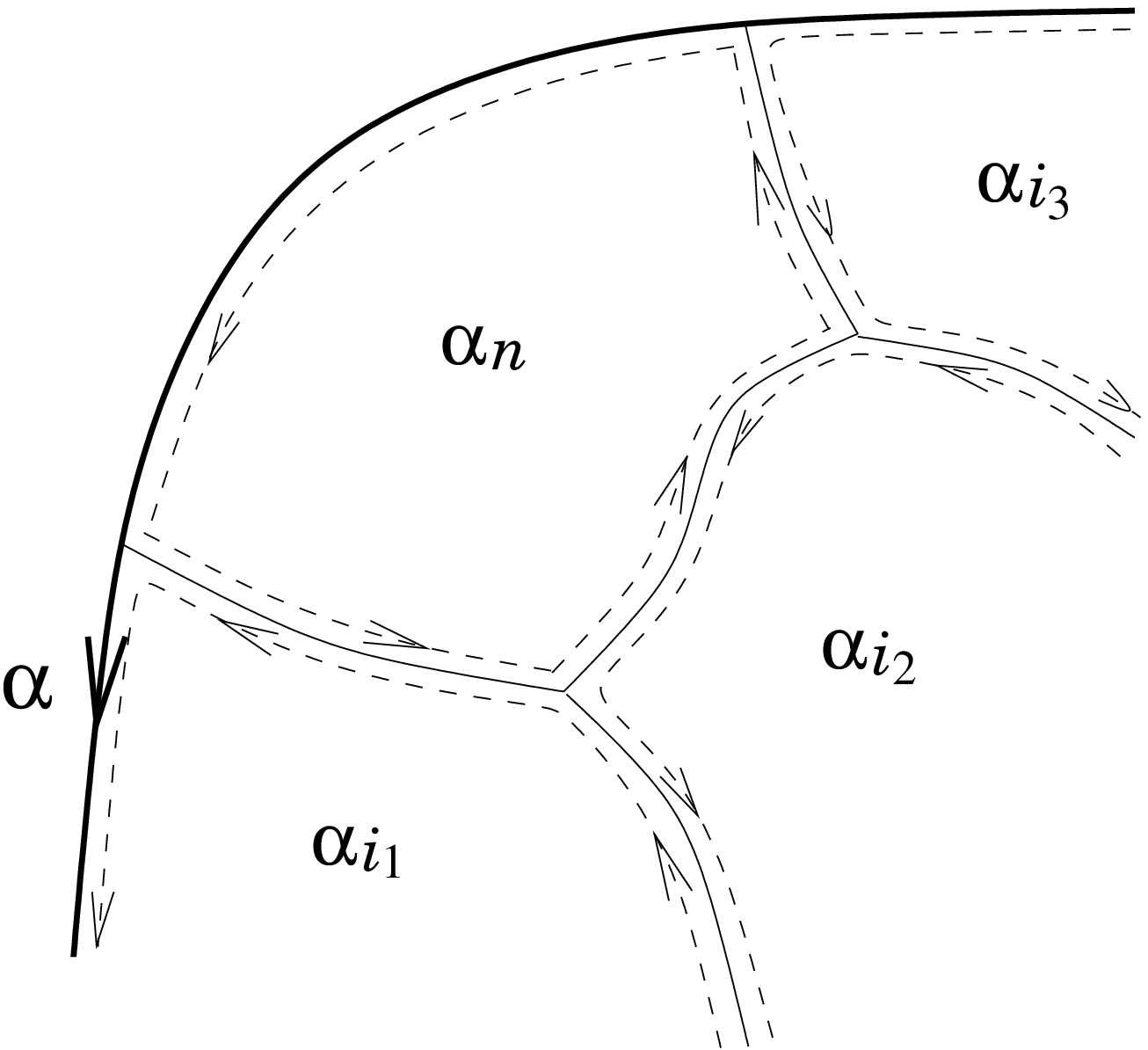}
\hbox to 5mm{}
\epsfxsize 55mm\epsfbox{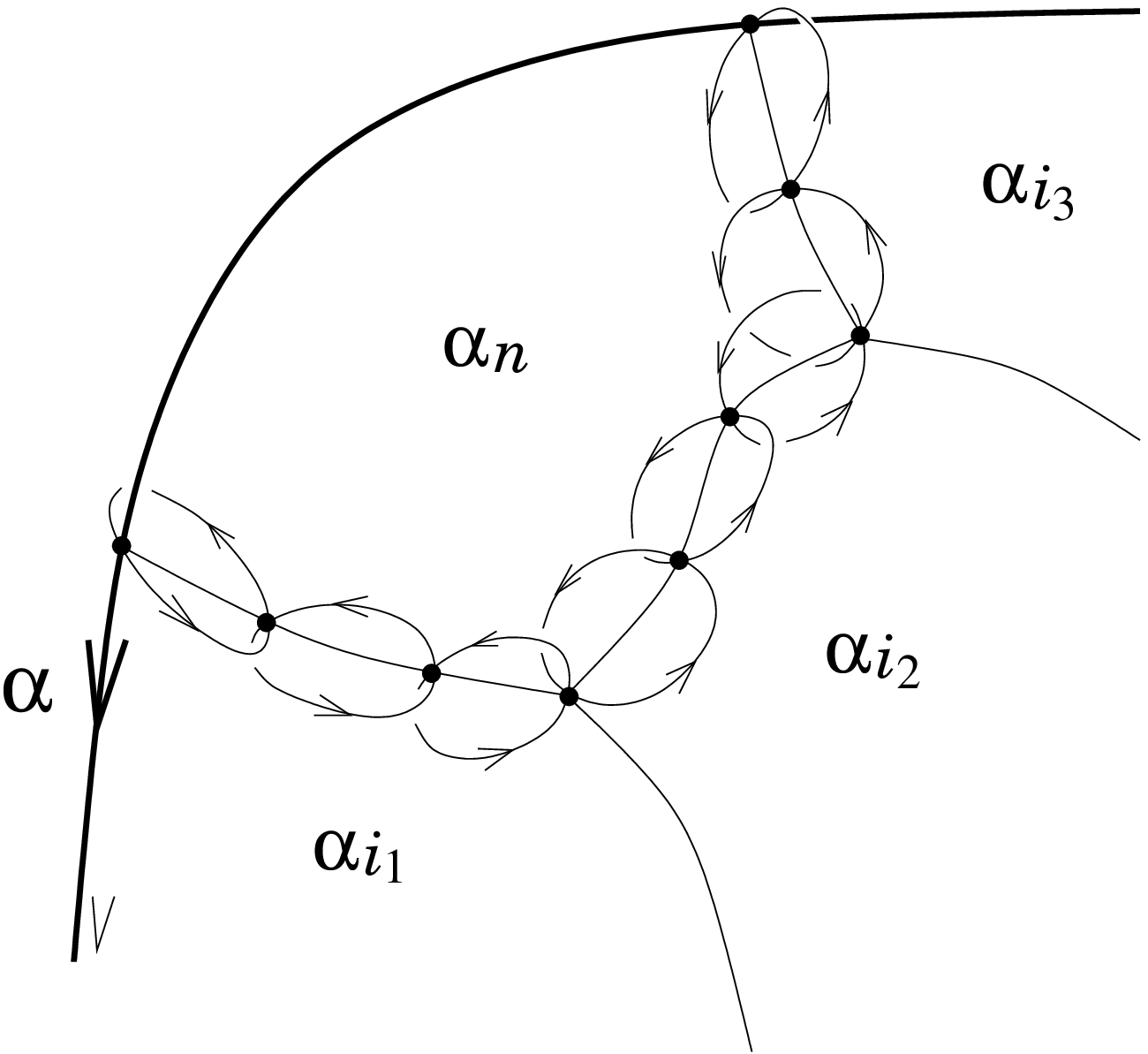}
}
\centerline{ 
\hbox to 20mm{}  
a). The net $\calA$.
\hbox to 25mm{}  
b). $(\supp\calA)\cup(S_1\cup\dots\cup S_1)$.
}
\medskip
\centerline{
\epsfxsize 55mm\epsfbox{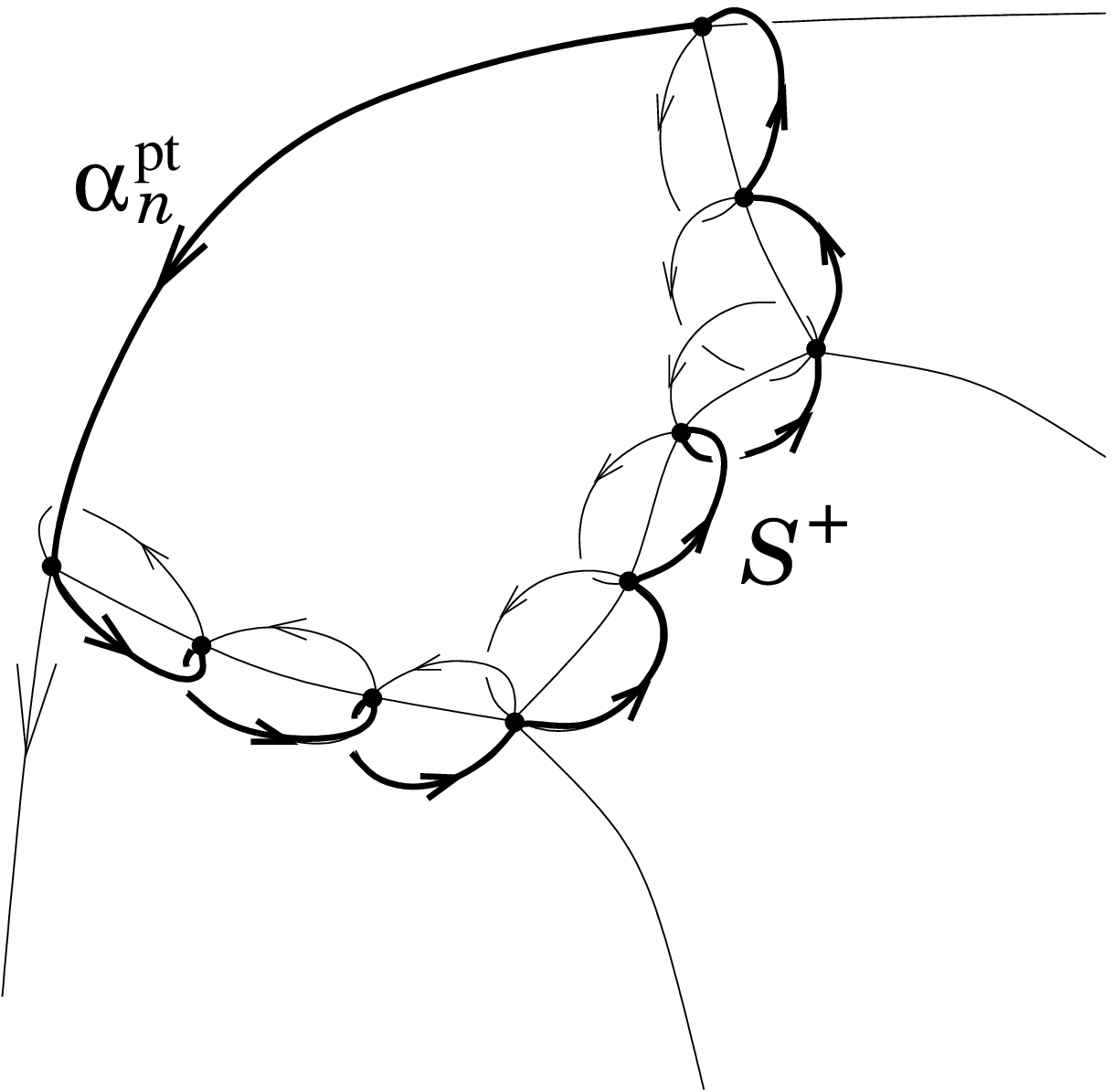}
\hbox to 5mm{}
\epsfxsize 55mm\epsfbox{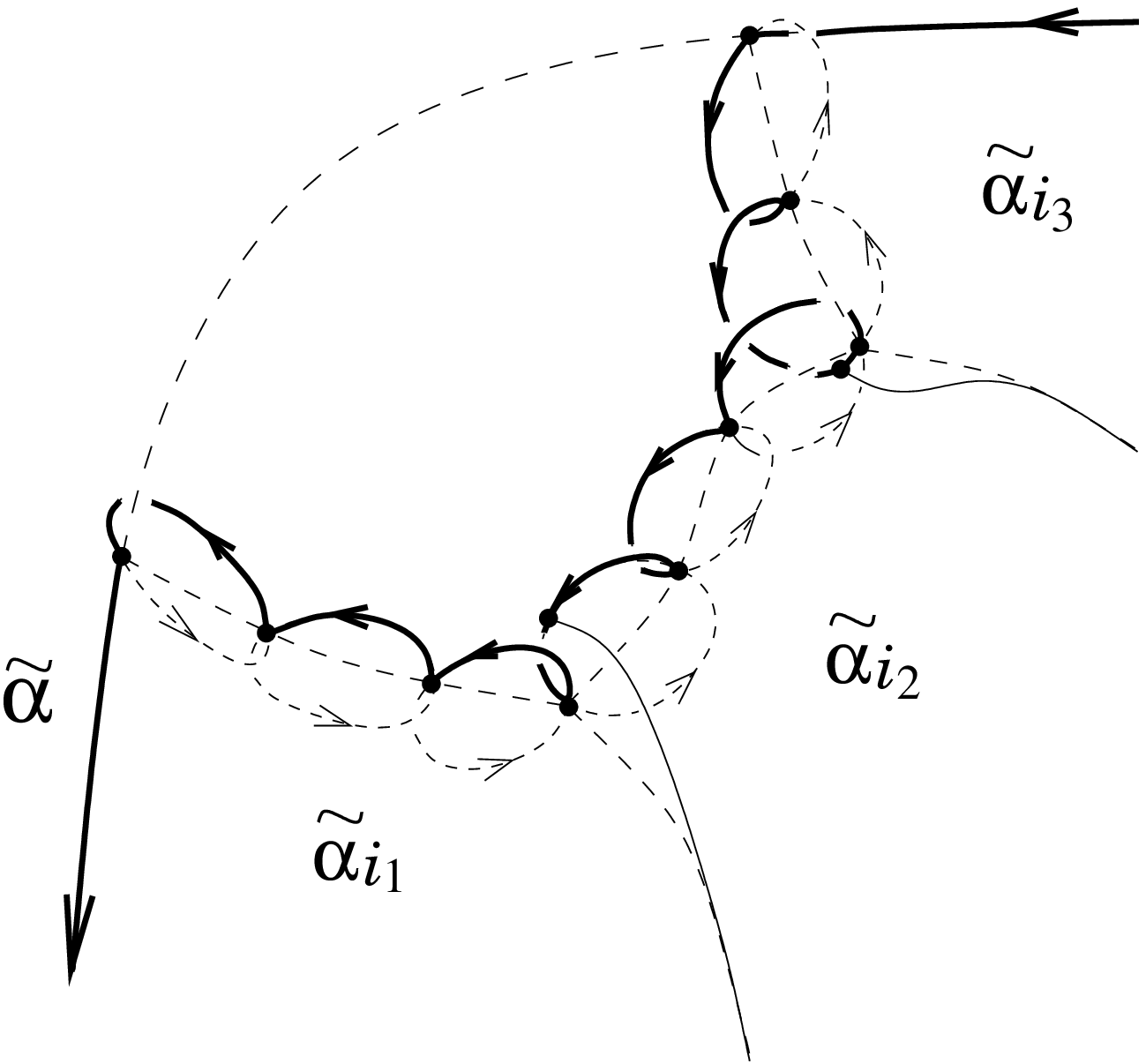}
}
\centerline{ 
c). The circuit $\beta_n = \alpha_n^{\pt} + S^+$.
\hbox to 30mm{}  
d). The net $\tilde\calA$.
\hbox to 10mm{}  
}
\botcaption{ Fig.~\figLegNet }
\endcaption
\endinsert

 By Corollary {\corAppr}, for any point
 $p\in\alpha_n^{\leg}$ there is a neighbourhood $U_p$ such that for any 
 $q\in U_p\cap\alpha_n$ there exists an analytic disk 
 $D_{pq}\subset\Omega$ satisfying the estimate (\eqApprSix) with an arbitrarily given
 number $\delta_1$ instead of $\delta$.
 Choosing a finite subcovering $\{U_p\}$, we can represent $\alpha_n^{\leg}$
 as the sum of arcs $\alpha_n^{\leg} = \gamma_1+\dots+\gamma_k$ so that for any  
 $i=1,\dots,k$ there exists an analytic disk
 $D_i\subset\Omega$ such that $\partial D_i = S_i = S_i^+ + S_i^-$, 
 $\partial S_i^\pm=\pm\partial\gamma_i$, and 
 $\length(S_i^\pm)/\length(\gamma_i) < 1+\delta_1$.
 We may also assume that the length of each arc $\gamma_i$ is less than an 
 arbitrarily given number, and that any edge of $\Gamma$
 (the graph from Definition \defNetwork)
 contributing to $\alpha_n^{\leg}$ is the sum of a subset of arcs $\gamma_i$.

 Perturbing the disks $D_i$, we can achieve that they are transverse to each other
 and hence, the curves $S_i$ have distinct tangents at the intersection points.
 We may also assume that if an end of 
 $\gamma_i$ lies on $\alpha$ then the tangents at this point to 
 $\alpha$ and to $\gamma_i$ are distinct.
 Let us set $S^\pm = \sum_{i=1}^k S_i^\pm$. 
 These are positively transverse chains such that 
 $\partial S^+=\partial\alpha_n^{\leg}=-\partial S^-$.
 Hence $\tilde\alpha = \alpha-\alpha_n^{\leg}+S^-$ is a generic positively
 transverse cycle.

 Passing if necessary from the arcs $\gamma_j$ to their subdivisions, we
 may assume that the collection of arcs $\gamma_1,\dots,\gamma_k$ can be
 completed up to 
 $\gamma_1,\dots,\gamma_k,\gamma_{k+1},\dots,\gamma_m$,
 so that $\alpha_i = \sum_{j=1}^m a_{ij} \gamma_j$, $i=1,\dots,n$,
 for some matrix of integer coefficients $a_{ij}$ such that 
 $a_{ij}\in\{-1,0,1\}$ for any $i=1,\dots,n$, $j=1,\dots,m$.
 Some of $\gamma_{k+1},\dots,\gamma_m$ being positively transverse, 
 the others being Legendrian.
 
Let us denote the set of ends of $\gamma_1,\dots,\gamma_k$ not belonging to
$\alpha$ by $P$. In other words, 
$P=(S\cap\supp\alpha_n)\setminus\supp\alpha=
(S^-\cap\supp\alpha_n)\setminus\supp\alpha$.
For every $p\in P$, let us define a point $\tilde p$ as follows.
Let $\gamma_i$, $1\le i\le k$, be an arc whose end is $p$
(there are two such arcs but we choose any of them).
Then we define $\tilde p$ as an interior point of $S^-_i$ which is closer to
$p$ than to the other end of $S^-_i$. 
If $p$ is the end of an arc $\gamma_i$ and $p\not\in P$,
we set $\tilde p = p$.

For any $i=1,\dots,m$, let us define an arc $\tilde\gamma_i$ as follows.
Let $\partial\gamma_i=q-p$, and let $\tilde p$ and $\tilde q$ be the points
chosen as it is described above starting from $p$ and $q$ respectively.
If $1\le i\le k$, we define $\tilde\gamma_i$ as the path on $S^-$ 
connecting $\tilde p$ to $\tilde q$.
If $k<i\le m$ and the arc $\gamma_i$ is Legendrian,
we define $\tilde\gamma_i$ as a Legendrian path from 
$\tilde p$ to $\tilde q$. If the arc $\gamma_i$ is positively transverse,
we set $\tilde\gamma_i=\gamma_i$.
In all the cases, we orient $\tilde\gamma_i$ so that 
$\partial\tilde\gamma_i = \tilde q-\tilde p$.
It follows from Lemma {\lemLegTwo} that the arc
$\tilde\gamma_i$ can be chosen arbitrarily close to $\gamma_i$.

Let us set $\tilde\calA = \{\tilde\alpha_1,\dots,\tilde\alpha_{n-1}\}$
where $\tilde\alpha_i = \sum_{j=1}^m a_{ij} \tilde\gamma_j$, $i=1,\dots,n-1$.
It is easy to check that this is a generic Legendrian net hanged on $\tilde\alpha$
(see Figure \figLegNet\,d).
Hence, by the induction hypothesis, we can find a PTSC-curve
$\tilde S = S_{k+1}+\dots+S_N$ so that the statement of the lemma holds for
$\tilde\calA$ instead of $\calA$ and for an arbitrarily chosen constant instead
of $\delta$. Then, for a suitable choice of the constants involved in the
construction of $\tilde\calA$, the curve 
$S = S_1+\dots+S_k + \tilde S$ will satisfy the conclusion of the lemma.
Indeed, let us denote the circuits of $\tilde S$ by 
$\beta_1,\dots,\beta_{n-1}$.
Then the curve $S$ has $n$ circuits, namely, $\beta_1,\dots,\beta_{n-1}$, and 
$\beta_n = \alpha_n^{\pt}+S^+$.
By the induction hypothesis, the circuits $\beta_1,\dots,\beta_{n-1}$
are close to the circles $\tilde\alpha_1,\dots,\tilde\alpha_{n-1}$,
hence also to the cycles $\alpha_1,\dots,\alpha_{n-1}$.
The circuit $\beta_n$ is close to the cycle $\alpha_n$ by construction
(see Figure \figLegNet\,c).

\medskip

b).
The proof if more or less the same as in Part a), but the manifold $M$ should
be replaced by a neighbourhood of the support of $\calA$ where the
quantity $K_1/K_2$
from (\eqAppr) is bounded from the below by some constant.
\qed\enddemo

Let us denote $\R_+=\{x\in\R\,|\,x\ge 0\}$ and $\R_-=\{x\in\R\,|\,x\le 0\}$.

\proclaim{ Lemma \lemPerturbLin } Let $x,y,z$ be coordinates in $\R^3$
and let $f_1,f_2:\R^3\to\C^2$ be the mappings given by 
$f_1(x,y,z) = x+iz$, $f_2(x,y,z) = y+iz$. For any complex number $c$,
let us denote the real curve $\{p\in\R^3\,|\,f_1(p)f_2(p)=c\}$ by $S_c$.
Then, for $c\not\in\R_-$, the curve $S_c$ has exactly two branches 
{\rm(}i.e. two connected components{\rm)}
$S_c^+$ and $S_c^-$ such that $S_c^+\subset\{x+y>0\}$, $S_c^-\subset\{x+y<0\}$,
and the restriction of the linear function $\R^3\to\R$, $(x,y,z)\mapsto x-y$, 
to each of the branches $S_c^\pm$ is a diffeomorphism.

Moreover, $S^\pm_c$ tends in any reasonable sense to $S^\pm_0$ as 
$c\to0$, $c\not\in\R_-$ where $S_0^+=\{z=xy=0, x+y\ge0\}$ and $S_0^-=\{z=xy,x+y\le0\}$.
\endproclaim

\demo{ Proof }
Set $a=\Re c$, $b=\Im c$. Then the curve $S_c$ is given by the system of
simultaneous equations $xy-z^2=a$, $z(x+y)=b$. By the change of variables 
$x-y=2u$, $x+y=2v$ we transform this system to $v^2-u^2-z^2=a$, $2zv=b$.

If $b=0$ and $a>0$ then $S_c$ is the hyperbola $v^2-u^2=a$ in the plane $z=0$.

If $b\ne 0$ then the intersection of $S_c$ with the plane $u=u_0$ can be found
by solving the system of simultaneous equations
$v^2-u^2-z^2=a$, $2zv=b$, $u=u_0$. Eliminating 
$u,z$, we obtain the equation $v^4 - (u_0^2+a)v^2 - (b/2)^2 = 0$ with respect to 
the variable $v$. It is clear that for any value of $u_0$, this equation has
exactly two roots one of whom being positive and the other one being negative. 
\qed\enddemo

\remark{ Remark }
For $c\in\R_-$, the curve $S_c$ is not smooth. It is the union of the hyperbola
$v^2-u^2=c$ in the plane $z=0$ and the circle $u^2+z^2=-c$ in the plane
$v=0$ which cross each other at the two points $z=v=0$, $u=\pm\sqrt{-c}$.
\endremark

\proclaim{ Lemma \lemPerturb }
Let $M$ be a $C^2$-smooth oriented real $3$-manifold, and let
$f_1,f_2,h$ be $C^2$-smooth complex valued functions on $M$ such that 
$f_1(p_0)=f_2(p_0)=0$, $h(p_0)\ne 0$, and
each of $f_1,f_2$ is a submersion at a neighbourhood of some point $p_0\in M$. 
Let us denote the real curves $f_j^{-1}(0)$ by $\gamma_j$, $j=1,2$. 
On each $\gamma_j$ near $p_0$, 
let us introduce the orientation induced by the submersion $f_j$.
Suppose that the tangents to $\gamma_1$ and $\gamma_2$ at $p_0$ are distinct.

Then there exist a number $\theta_0\in[0,2\pi[$ and a neighbourhood $U$ of $p_0$ such that
each of the curves $U\cap\gamma_j$, $j=1,2$, is diffeomorphic to an open interval
and for any fixed $\theta\not\equiv\theta_0\mod2\pi$ 
there exists $r_0=r_0(\theta)>0$ such that for $0<r<r_0$, the curve 
$S_{r,\theta} = \{\,p\in U\,|\, f_1(p)f_2(p)=re^{i\theta}h(p)\,\}$ consists of 
two smooth branches one of which tending to $\gamma_1^-\cup\gamma_2^+$ and the other
one tending to $\gamma_2^-\cup\gamma_1^+$ as
$r\to0$ where $\gamma_j^\pm$ denotes the preimage of $\R_\pm$ under an orientation
preserving embedding $(U\cap\gamma_j,p_0)\to(\R,0)$.
\endproclaim

\demo{ Proof } It is clear that if the statement of the lemma holds for
$f_1,f_2,h$, then it holds 
(with maybe another number $\theta_0$) also for
$c_1f_1,c_2f_2,c_3h$ where $c_1,c_2,c_3$ are arbitrary nonzero complex numbers.
Therefore, we may assume that $h(p_0)=1$. Let us choose a local real coordinate
$z$ in a neighbourhood of $p_0$ so that the both curves $\gamma_1$, $\gamma_2$ 
lye on the surface $z=0$. Multiplying $f_1$ and $f_2$ by suitable complex numbers,
we may assume that $\partial(\Re f_j)/\partial z(p_0)=0$, $j=1,2$.
Let us set $x=\Re f_1$, $y=\Re f_2$. Then $(x,y,z)$ is a local coordinate system 
where the functions $f_1,f_2$ have the form
$f_1(x,y,z) = x+iz + O(x^2+y^2+z^2)$, $f_2(x,y,z) = y+iz + O(x^2+y^2+z^2)$.
Therefore, the statement follows from Lemma {\lemPerturbLin} combined with the fact that
the curve $H_r( S_{r,\theta})$ tends to the curve $\{(x+iz)(y+iz)=e^{i\theta}\}$ as 
$r\to 0$, $\theta=\const\not\equiv\pi\mod2\pi$ where $H_r$ stands for the 
homothety
$(x,y,z)\mapsto(x/\sqrt r,\,y/\sqrt r,\,z/\sqrt r)$.
\qed\enddemo

\proclaim{ Proposition \propPerturb }
Let $\Omega$ be an arbitrary domain in $\C^2$ with a compact $C^2$-smooth boundary
$M=\partial\Omega$, and let $A$ be an algebraic curve in $\C^2$ given by $f=0$.
Suppose that $S=A\cap M$ is a PTSC-curve. Let $h$ be a polynomial which does not vanish
at the double points of $S$. Then there exists a finite set
$\Theta\subset[0,2\pi[$ such that for any  
$\theta\in[0,2\pi[\,\setminus\Theta$, there is $r_0=r_0(\theta)>0$ such that
the real curve 
$S_{r,\theta} = \{\,p\in M\,|\, f(p)=re^{i\theta}h(p)\,\}$ for $0<r<r_0$ is smooth
and its connected components converge to the circuits of $S$ as $r\to 0$.
\endproclaim

\demo{ Proof }
Follows from Lemma \lemPerturb.
\qed\enddemo

 \demo{ Proof of Theorem \thMain }
By Proposition {\propLeg}, we can construct a Legendrian net hanged on 
$\partial A$ all whose cells are small. 
By Proposition {\propAppr}, it can be approximated by the boundary of the union $D$
of analytic (resp. linear) disks so that the circuits of $D\cap M$ are arbitrarily
small. Using Proposition {\propPerturb} in the algebraic case, and the standard
techniques of analytic sheaves on open Riemann surfaces in the analytic case,
we can perturb $D$ so that all its boundary components become close to circuits
of $S$. Moreover, this perturbation can be chosen so that the points of $P$
do not move (in the algebraic case, we just choose $h$ in Proposition 
{\propPerturb} which does not vanish on $P$).
 \qed
 \enddemo

%
%
%

\head \sectStokes. Some elementary facts about the standard contact structure
                   on $\SS^3\subset\C^2$.
\endhead

For the reader's convenience, in this section we shall give some well-known facts
about curves on $\SS^3$ and their projections to $\P^2$, and we shall deduce
from them a lower bound for  $n(\eps)$ of the order $1/\eps^2$.

Let $z=x+iy$ and $w=u+iv$ be the standard coordinates in $\C^2$.
Let us denote:
$$
          \rho = \rho(z,w) = |z|^2 + |w|^2 = x^2+y^2+u^2+v^2,
$$
$$
        \eta = 1/2\,d^c\rho=i/2\,(z\,d\bar z - \bar z\,dz + w\,d\bar w - \bar w\,dw)
        = x\,dy - y\,dx + u\,dv - v\,du,
$$
$$
        \omega = 1/2\,d\eta = i/2\,(dz\wedge d\bar z + dw\wedge d\bar w)
        = dx\wedge dy + du\wedge dv.
$$
Let
$\B^4=\{\rho\le 1\}$, $\SS^3 = \partial\B^4 = \{\rho=1\}$,
$\P^1 = \big(\C^2\setminus\{(0,0)\}\big)/_{(z,w)\sim(\lambda z,\lambda w)}\,$  
and let $\pr:\C^2\setminus\{(0,0)\}\to\P^1$ 
and $\pr_{\SS^3}:\C^2\setminus\{(0,0)\}\to\SS^3$ be the standard projections.
The field of real 2-planes $\ker\eta|_{\SS^3}$ is the field of complex tangents to 
$\SS^3$. It defines the standard (tight) complex structure on $\SS^3$.

Let $\|\cdot\|_{\P^1}$ and $\omega_{\P^1}$ be the Riemannian Fubini-Studi metric on
$\P^1$ and the corresponding volume form which are defined by
$$
        \|d\zeta\|_{\P^1} = {|d\zeta|^2\over(1+|\zeta|^2)^2},\qquad
        \omega_{\P^1} = {i\over2}\,{d\zeta\wedge d\bar\zeta\over(1+|\zeta|^2)^2},
        \qquad  \zeta = z/w.
$$
$\P^1$ equipped with this metric is isometric to the standard $2$-sphere of the radius 
$1/2$, in particular, we have
$$
 \int_{\P^1} \omega_{\P^1} = \pi.
$$

Let
$$
     \eta^* = \pr_{\SS^3}^*(\eta|_{\SS^3})
     \qquad\text{and}\qquad
        \omega^* = \pr^*(\omega_{\P^1}).
$$
It is easy to check that  
$$
 \eta^* = {\eta\over\rho} = {1\over2}\,d^c\log\rho
     \qquad\text{and}\qquad
 d\eta^* = {2\omega\over\rho} - {d\rho\wedge\eta\over\rho^2} = 2\omega^*.
                                                                        \eqno(\eqdrho)
$$

\proclaim{ Lemma \lemStokes }
Let $F$ be a $2$-chain in $\SS^3$. Then 
$$
      \int_{\partial F} \eta = 2\int_{\pr_*F} \omega_{\P^1}.
$$
\endproclaim

\demo{ Proof } Follows from Stokes' theorem and from (\eqdrho).
\qed\enddemo

Let $\|\cdot\|_{\SS^3}$ be the Riemannian metric on $\SS^3$ induced by the standard
metric in $\C^2$. It is easy to check that 
$$
     \|v\|_{\SS^3}^2 = |\eta(v)|^2 + \|\pr_* v\|_{\P^1}^2,\qquad v\in T\SS^3.
                                                                    \eqno(\eqNorm)
$$
In particular, if $D$ is the disk which is cut on $\B^4$ by a complex line
passing through the origin, then the circle  $\partial D$ is orthogonal to the
contact structure and
$$
      \int_{\partial D}\eta = 2\pi.                          \eqno(\eqTwoPi)
$$

\proclaim{ Lemma \lemDA }
Let $A$ be a smooth complex algebraic curve in $\C^2$ passing through the origin
and having there a non-degenerate tangency with a complex line $L$.
Let $F$ be the closure of  
$\pr_{\SS^3}(A\cap\B^3\setminus 0)$.
Then $\partial F = \partial(A\cap\B^4) - \partial(L\cap\B^4)$.
In particular, 
$$
    \int_{\partial(A\cap\B^4)}\eta = 2\pi + 2\int_{\pr_* F}\omega_{\P^2}
     \ge 2\pi.
                                                             \eqno(\eqStokes)
$$
\endproclaim

\demo{ Proof }
Apply the {\sl real} blowup of the origin
(identifying $\C^2$ with $\R^4$). 
\qed\enddemo


\definition{ Definition \defBetaPN }
An $n$-chain $\beta$ with a piecewise smooth boundary on an oriented
$n$-manifold $M$ is called {\it positive} (resp. {\it strictly positive})
if each connected component of the complement of 
$\supp\beta$ contributes to $\beta$ 
with a {\it nonnegative} (resp. {\it positive}) multiplicity.
We shall write in this case $\beta\ge0$ (resp. $\beta>0$).

Every $n$-chain $\beta$ on $M$ can be represented in a unique way as
$\beta = \beta^+ - \beta^-$ so that $\beta^+\ge0$, $\beta^-\ge0$,
and $(\supp\beta^+)\cap(\supp\beta^-)=(\supp\partial\beta^+)\cap
(\supp\partial\beta^-)$. The chains $\beta^\pm$ are called the  
{\it positive} and the {\it negative parts} of $\beta$.

If $U$ is a domain in $M$ which has a piecewise smooth boundary and if
$\beta$ is an $n$-chain then the {\it restriction} of $\beta$ to $U$ is
the $n$-chain $\beta|_U = \sum m_i(\beta_i\cap U)$ where 
$\beta = \sum m_i\beta_i$ is the representation of
$\beta$ as a linear combination of domains with piecewise smooth boundaries.
в виде целочисленной линейной комбинации областей
\enddefinition

\remark{ Remark \remBetaPN }
Let $M$ be an oriented $n$-manifold.
We shall identify $n$-chains on $M$ having piecewise smooth boundaries
with integer-valued functions being linear combinations of characteristic
functions of domains. Namely, if  $\beta_1,\dots,\beta_k$ are domains in  
$M$ having piecewise smooth boundaries then the chain
$\beta = \sum m_i\beta_i$, $m_i\in\Z$, will be identified with the function 
$\chi_\beta = \sum m_i\chi_{\beta_i}$ where $\chi_{\beta_i}$ is the 
characteristic function of the domain $\beta_i$ (i.e. $\chi_{\beta_i}|_{\beta_i}=1$,
$\chi_{\beta_i}|_{M\setminus\beta_i} = 0$).

The integral of a 2-form $\xi$ corresponds under this identification to
$\int_M\chi_\beta\xi$. Taking the restriction of $\beta$ to $U$ corresponds to
the multiplication by  $\chi_U$, etc.
\endremark

\proclaim{ Lemma \lemIsop } 
{\rm(}Isoperimetric inequality for $2$-chains on $\SS^2$.{\rm)}
Let $\SS^2$ be the sphere of a radius $R$ in $\R^3$ endowed with the standard
Riemannian metric and the standard area form $dS$.
Let $\beta$ be a $2$-chain on $\SS^2$ which has a piecewise smooth boundary
whose length  {\rm(}taking into account the multiplicities if there are multiple
segments{\rm)} is equal to $a$, and let
$b=\int_\beta dS$ be the oriented area of $\beta$.
Let $\beta^+$ {\rm(}resp. $\beta^-${\rm)} be the positive
 {\rm(}resp. negative{\rm)} part of $\beta$, and let
$b^\pm=\int_{\beta^\pm}dS$.

Suppose that $|b| < 2\pi R^2$ и $a < 2\pi R$. Then
$$
   |b| \le b^+ + b^- \le S_R(a) \qquad\text{where}\quad
   S_R(a) = 2\pi R^2\big(1-\sqrt{1-{a^2/ (2\pi R)^2}}\,\big)  \eqno(\eqIsop)
$$
and if, moreover, the set $\supp\partial\beta$ is connected then
$$
   \diam_{\SS^2}\supp\beta\le a/2.                                    \eqno(\eqDiam)
$$
\endproclaim

\demo{ Proof }
If $\beta$ is a domain on the sphere then (\eqIsop) is the classical 
isoperimetric inequality.

In the general case, the boundary of $\beta$ can be represented as a disjoint union
of closed curves whose lengths we denote by $a_1,\dots,a_k$.
Each of these curves is the common boundary of two domains in the sphere,
and the area of at least one of them does not exceed 
$2\pi R^2$. Choosing in a suitable way the signs of these domains, we obtain
a 2-chain whose boundary coincides with $\partial\beta$.
Adding if necessary several times $\pm[\SS^2]$, we obtain a 2-chain
$\beta'$ such that $\beta-\beta'$ is a zero-homologous cycle,
$\partial\beta'=\partial\beta$, and $\beta'$ has the form
$\beta'=m[\SS^2]+s_1\beta_1+\dots+s_k\beta_k$ where $m\in\Z$,
$s_i=\pm1$, and
$\beta_i$ is a domain of area $b_i\le 2\pi R^2$.
For each of these domains, we have $b_i \le S_R(a_i)$.
Since the function $S_R$ is convex and $S_R(0)=0$,
it follows that $S_R(a)=S_R(a_1+\dots+a_k)\ge S_R(a_1)+\dots+S_R(a_k)$.
Hence,
$$
  |b-4\pi R^2m| = \Big|\sum s_i b_i\Big| \le \sum b_i \le \sum S_R(a_i)
   \le S_R(a).
$$
Let us show that $m=0$. Indeed, recall that
$|b|<2\pi R^2$. Combining this inequality with $S_R(a) < 2\pi R^2$, we obtain
$4\pi R^2 |m| \le |b| + |b-4m\pi R^2| < 2\pi R^2 + S_R(a) < 4\pi R^2$, i.e. $|m|<1$.
But $m\in\Z$, hence $m=0$.

Let us set $\hat\beta^\pm = \sum_{s_i=\pm1}\beta_i$, $\hat b^\pm = \sum_{s_i=\pm1}b_i$.
We have proven that $\hat b^+ + \hat b^- \le S_R(a)$ and for deducing
(\eqIsop), it remains to note that $b^\pm\le\hat b^\pm$. The latter fact is evident
because the decomposition  
$\beta^+-\beta^-$ can be obtained from $\hat\beta^+ - \hat\beta^-$ by a successive 
cancellation of connected components of $\supp\partial\beta$ contributing
simultaneously to $\beta^+$ and $\beta^-$.

Now let us suppose that the set $\supp\partial\beta$ is connected and let us prove
(\eqDiam).
First, let us show that $\supp\beta$ does not contain any pair of antipodal points.
Indeed, let us denote the central symmetry by 
$\sigma:\SS^2\to\SS^2$. 
The estimate $\length\partial\beta<2\pi R$ yields
$\sigma(\supp\partial\beta)\cap\supp\partial\beta=\varnothing$.
Since $\supp\beta$ is connected, this implies that 
$\sigma(\supp\beta)$ is contained in a single connected component of
the complement of $\supp\partial\beta$. This component cannot be contained in
$\supp\beta$ because its area is greater than the area of
$\sigma(\supp\beta)$, hence, it is greater than the area of $\supp\beta$.
Therefore, we have $\sigma(\supp\beta)\cap \supp\beta=\varnothing$.

Let $p,q\in\supp\beta$. Let us denote the shortest geodesic from $p$ to $\sigma(q)$
(resp. from $q$ to $\sigma(p)$) by $\gamma_p$ (resp. by $\gamma_q$).
Since the points $\sigma(p)$ and $\sigma(q)$
do not belong to $\supp\beta$, there exist points 
$p'\in\gamma_p\cap\supp\partial\beta$ and
$q'\in\gamma_q\cap\supp\partial\beta$.
Therefore, we have
$\dist_{\SS^2}(p,q)\le\dist_{\SS^2}(p',q')\le (\length\partial\beta)/2=a/2$.
\qed\enddemo

\remark{ Remark \remIsop } The classical isoperimetric inequality
(the inequality (\eqIsop) for a single domain in the sphere)
can be equivalently reformulated as
$4\pi b - b^2/R^2 \le a^2$. In this form, it holds without the assumptions
$b < 2\pi R^2$ and $a < 2\pi R$.
An analogue of this inequality for $2$-chains is
$$
   \max_{m\in\Z}\Big( |b-4m\pi R^2| - (b-4m\pi R^2)/R^2\Big) \le a^2.
$$
It holds also without the assumptions
$|b| < 2\pi R^2$ and $a < 2\pi R$.
The graph of the left hand side of the latter inequality (considered as a
function of $b$) is the union of the upper halves of ellipses centered at
the points
$\big((2+4m)\pi R^2,0\big)$, $m\in\Z$. The ellipses touch each other at the points
$\big(4m\pi R^2,0\big)$.
\endremark

\proclaim{ Lemma \lemLength }
Let $\gamma$ be a positively transverse curve on $\SS^3$
{\rm(}e.g. a connected component of the intersection of a complex analytic curve with
$\SS^3${\rm)}.
Let us denote:
$$
    a = \length_{\P^1}(\pr\gamma),\qquad
    b = \int_\gamma\eta, \qquad
    \ell = \length_{\SS^2}(\gamma).
$$
Then we have
$$
   \max(a,b) \le \ell \le a+b                         \eqno(\eqABL)
$$
and if, moreover, $\ell<\pi/2$ then we have
$$
    b \le S_{1/2}(a) = 
    {\pi\over 2}\,\Big( 1 - \sqrt{ 1 - {a^2\over\pi^2}\! }\,\,\Big)
    = {a^2\over 4\pi} + {a^4\over 16\pi^3} +\dots            \eqno(\eqMain)
$$
\endproclaim

\demo{ Proof } The inequalities (\eqABL)
follow from (\eqNorm) combined with the fact that the form 
$\eta$ is positive on $\gamma$.

To prove (\eqMain), let us consider a $2$-chain in 
$\SS^3$ whose boundary is the cycle $\gamma$ and let us denote its projection 
to $\P^1$ by $\beta$. Recall that $\P^1$ is isometric to the sphere of the radius 
$R=1/2$.
Hence, $\ell<\pi/\sqrt2$ combined with (\eqABL) implies
$a < \ell < \pi/2 < \pi = 2\pi R$ and 
$b < \ell < \pi/2 = 2\pi R^2$, and the result follows from Lemma {\lemIsop}.
\qed\enddemo

\proclaim{ Corollary \corLength } Let 
$\gamma$ be a positively transverse curve in $\SS^3$ and let 
$\ell$ and $b$ be as in Lemma \lemLength.
If $\ell < \pi/2$ then $b \le S_{1/2}(\ell)$.
\endproclaim

\demo{ Proof }
Follows from (\eqABL), (\eqMain), and the monotonicity of $S_{1/2}$
\qed
\enddemo

Combining all the above facts, we easily obtain the quadratic estimate for
$n(\eps)$ which we announced in Introduction:

\proclaim{ Proposition \propEpsSq } If $\eps < \pi/2$ then
$n(\eps) > 2\pi/S_{1/2}(\eps) = 8\pi^2/\eps^2 - 2 + O(\eps^2)$.
\endproclaim

\demo{ Proof }
Let $A$ be a complex algebraic curve in $\C^2$ passing through the origin
such that all connected components 
$\gamma_1,\dots,\gamma_n$ of $A\cap\SS^3$ are shorter than $\eps$.
Perturbing $A$, we may assume that the conditions of Lemma {\lemDA} are satisfied.
Thus, by (\eqStokes) and Corollary \corLength, we have
$$
   2\pi\le\int_{\partial(A\cap\B^4)}\eta = \sum_{i=1}^n\int_{\gamma_i}\eta
   \le n S_{1/2}(\eps) = 
   {n\eps^2\over 4\pi}\cdot\big(1+{\eps^2\over4\pi^2}+O(\eps^4)\big). \qed
$$
\enddemo

%
%

\head \sectLegLB. A lower bound of the order $\eps^{-3}$ for the number of
                  cells of a Legendrian net
\endhead

In this section, we shall prove the following result which means 
that any upper bound obtained by the method of {\S\S\sectLeg--\sectAppr}
cannot be better than $n(\eps) = O(\eps^{-3})$.
More precisely, we shall prove the following result.

\proclaim{ Proposition \propLegLB } Let $A$ be an algebraic curve in 
$\C^2$ passing through the origin, and let $\Gamma = A\cap\SS^3$.
Let $\calA=\{\alpha_1,\dots,\alpha_n\}$ be a Legendrian net hanged on
$\Gamma$ {\rm(}see Definition
\defNetwork{\rm)}. Suppose that every cell of $\calA$ is shorter than $\eps$.
Then 
$$
     n > {2 c_0 \over \eps S_{1/2}(\eps)} = {8c_0\pi\over\eps^3} 
       - {2c_0\over\pi\eps} + O(\eps)
$$
where $c_0$ is a constant depending only on $A$. In the case when
$A$ is a complex line, one can set $c_0 = \pi^2/4$.
\endproclaim

\remark{ Remark } It seems that a similar statement should take place for any
contact $3$-manifold.
\endremark

\demo{ Proof }
We shall use the notation introduced in  \S\sectStokes.
We shall assume that $\eps<\pi/2$.
Let us define a function $f:\P^1\to\R_+$ by setting $f(p) = \dist_{\P^1}(p,\pr\Gamma)$.
For each cell $\alpha_i$, let us consider a 2-chain  $\tilde\beta_i$
in $\SS^3$ such that $\partial\tilde\beta_i=\alpha_i$, and let us set
$\beta_i = \pr_*\tilde\beta_i$. Let $\beta_i=\beta_i^+ - \beta_i^-$ be the decomposition
of  $\beta_i$ into the positive and the negative part
(see Definition \defBetaPN). Let us denote
$$
    b_i=\int_{\beta_i}\omega_{\P^1}, \quad
    b_i^+=\int_{\beta_i^+}\omega_{\P^1}, \quad
    b_i^-=\int_{\beta_i^-}\omega_{\P^1},\qquad
    c_0 = \int_{\P^1} f\,\omega_{\P^1}.
$$
In the case when $A$ is a complex line, it is not difficult to compute in
spherical coordinates that $c_0 = \pi^2/4$.
It follows from (\eqMain) that $b<\pi/2$, hence, by Lemma
{\lemIsop} we have
$$
    \diam_{\P^1}(\supp\beta_i)<\eps/2,\qquad  i=1,\dots,n.      \eqno(\eqDiamLeg)
$$
Perturbing $A$ if necessary, we may assume that it is non-degenerate.
Let $F$ and $L$ be as in Lemma \lemDA.
Let $\tilde\beta_0$ be a 2-chain in $\SS^3$ such that
$\partial\tilde\beta_0 = \partial(L\cap\B^4)$.
Then, according to (\eqStokes), we have
$$
    \sum_{i=1}^n \partial\tilde\beta_i = \sum_{i=1}^n \alpha_i = \Gamma = 
    \partial F + \partial(L\cap\B^4) = \partial F+\partial\tilde\beta_0.
$$
It follows from (\eqTwoPi) and from $\pr_*\partial\tilde\beta_0=0$, that
$\pr_*\tilde\beta_0 = [\P^1]$.
Hence $\sum_{i=1}^n\beta_i = \pr_* F + \pr_*\tilde\beta_0 = \pr_* F + [\P^1]$.
Thus,
$$
   c_0 = \int_{\P^1} f\,\omega_{\P^1}
       \le \int_{\pr_*F} f\,\omega_{\P^1} + \int_{\P^1} f\,\omega_{\P^1}
       = \sum_{i=1}^n \int_{\beta_i} f\,\omega_{\P^1}.
                                                                \eqno(\eqLegLBOne)
$$
Let us denote $m_i^+ = \max_{\supp\beta_i^+} f$, $m_i^-=\min_{\supp\beta_i^-} f$.
Then 

\noindent
\vbox{\vbox{
$$
   \int_{\beta_i} f\,\omega_{\P^1} = 
   \int_{\beta_i^+} f\,\omega_{\P^1} - 
   \int_{\beta_i^-} f\,\omega_{\P^1} \le 
   b_i^+ m_i^+ - b_i^- m_i^-   \qquad\qquad\qquad\qquad
$$}
\vbox{
$$
   \qquad\qquad\qquad =b_i^+(m_i^+ - m_i^-)
     + (b_i^+ - b_i^-) m_i^- 
   = b_i^+(m_i^+ - m_i^-)
     + b_i m_i^-.                                      \eqno(\eqLegLBTwo)
$$}}

\noindent
By Lemma {\lemIsop}, we have $b_i^+\le S_{1/2}(\eps)$.
Since $|f(p)-f(q)|\le\dist_{\P^1}(p,q)$, it follows from (\eqDiamLeg) 
that
$m_i^+ - m_i^- < \diam_{\P^1}\supp\beta_i < \eps/2$,
hence
$$
        b_i^+(m_i^+ - m_i^-) \le S_{1/2}(\eps)\eps/2.
                                                         \eqno(\eqLegLBThree)
$$
Let us show that 
$$
       b_i m_i^- = 0.
                                                         \eqno(\eqLegLBFour)
$$
Indeed, let $\alpha_i = \alpha_i^{\leg} + \alpha_i^{\pt}$ be the decomposition
from Definition \defNetwork. Let us consider two cases: $\alpha_i^{\pt}=0$ and
$\alpha_i^{\pt}\ne0$. 
In the former case, the cycle $\alpha_i$ is Legendrian, hence
$$
   b_i = \int_{\beta_i}\omega_{\P^1} = \int_{\tilde\beta_i}\omega^* 
       = \int_{\alpha_i}\eta^* = 0.
$$
In the latter case, $\supp\tilde\beta_i$ has a non-empty intersection with $\Gamma$,
hence $f$ vanishes in $\supp\beta_i$ which implies
$m_i^- = 0$. The equality (\eqLegLBFour) is proved.
Combining (\eqLegLBOne) -- (\eqLegLBFour), we obtain
$c_0 \le n S_{1/2}(\eps)\eps/2 = O(\eps^3)$.
\qed\enddemo

\remark{ Remark } In the case when $A$ is a complex line passing through the origin,
the quantity
$\int_{\beta_i}f\omega_{\P^1}$ playing the central role in the proof, can be 
interpreted as the moment of $\beta_i$ (considered as a measure on
$\P^1$) with respect to the point $\pr A$. 
So, the proof reduces to the following argument:
the measure $\omega_{\P^1}$ whose moment is equal to an absolute constant
$\pi^2/4$ is represented as the sum of measures
$\beta_i$ whose moments are of the order $\eps^3$.
\endremark

\medskip
Finally, let us formulate an open question, an affirmative answer to which
would imply a lower bound for  $n(\eps)$ of the order $\eps^{-3}$
(by the same method as Proposition {\propLegLB} is proved).

Let ${\Cal L}$ be the set of positive functions on 
$\P^1$ satisfying the Lipschitz condition with the constant 1, i.e.
functions  $f:\P^1\to\R_+$ such that
$|f(p)-p(q)| \le \dist_{\P^1}(p,q)$ for all $p,q\in\P^1$.

\smallskip
Does there exist an absolute constant  $c$ such that the inequality
$$   
    \max_{f\in {\Cal L}}\Big(
      \int_{m[\P^1]+\pr_*F} f\omega_{\P^1}
      - \int_{A\cap\SS^3} (f\circ\pr)\cdot\eta\Big) > c,
    \qquad
     F = \overline{ \pr_{\SS^3}(A\cap\B^4\setminus\{0\}) },
$$
holds for any algebraic curve $A\subset\C^2$ whose multiplicity at the origin
is $m$?
(As in (\eqStokes) and (\eqLegLBOne), here $\pr_*$ denotes the homomorphism
between the groups of 2-chains induced by $\pr:\SS^3\to\P^1$;
under the identification of 2-chains in $\P^1$ with integer-valued functions discussed
in Remark \remBetaPN, the 2-chain $m[\P^1]+\pr_*F$ corresponds to the function
whose value at the line $L$ through $0$, is the number of intersection points
of $L$ and $A\cap\B^4$ counting the multiplicities).


\head    \sectExplicit. Construction of a Legendrian net in 
         $\SS^3$ providing an upper bound for $n(\eps)$ of the order $1/\eps^3$
\endhead

Let us denote the coordinate axis $\{w=0\}$ by $L$. Let $\Gamma=L\cap\SS^3$.
For an integer $n$, we denote the rotation $(z,w)\mapsto(e^{2\pi i/n}z,w)$
by $\tilde R_n:\C^2\to\C^2$, and let 
$R_n:\P^1\to\P^1$ be the correspondent rotation 
$(z:w)\mapsto(e^{2\pi i/n}z:w)=(z:e^{-2\pi i/n}w)$.
Let us set $p_0 = (0:1)$, $p_\infty = (1:0)$. These are the fixed points of $R_n$.

Let us fix a small number  $\eps>0$, and let
$m=[10\pi/\eps]+1$.
Let us set 
$$
     r_k={k\pi\over 2m},\qquad
     \Delta_k = \{\,q\in\P^1\,|\,\dist_{\P^1}(p_0,q) \le 
     r_k\,\},
     \qquad
     k=0,\dots,m.
$$
Recall that $\P^1$ is isometric to the sphere of the radius $1/2$, hence
$$
   \{p_0\}=\Delta_0\subset\Delta_1\subset\dots\subset\Delta_{m} = \P^1.
$$
Let us denote the closure of $\Delta_k\setminus\Delta_{k-1}$  by $A_k$, and
let us set $a_k = \Area(A_k)$, $s_k = \Area(\Delta_k)=a_1+\dots+a_k$,
$k=1,\dots,m$. 
Let $\ell_k =(\length\partial\Delta_k+\length\partial\Delta_{k-1})/2$.
For each $k=1,\dots,m$, we set $n_k=2^{\nu_k}$ where $\nu_k$ is chosen so that
$$
     {\eps \over 40} < \ell_k^+  \le {\eps\over 20},
     \qquad 
     \ell_k^+={s_{k}\ell_k\over n_k a_k}.
                                                                      \eqno(\eqEll)
$$
It is clear that $\nu_k$ is uniquely determined by this condition.
Indeed, (\eqEll) implies that 
$\nu_k = [\log_2(20s_k\ell_k/(\eps a_k))]$.

By definition, we have
$$
     \ell_k={\pi\over2}\Big(\sin{k\pi\over m} + \sin{(k-1)\pi\over m}\Big),
     \quad
      s_k = {\pi\over2}\Big(1-\cos{k\pi\over m}\Big),
     \quad
      a_k = s_k - s_{k-1}.
                                                           \eqno(\eqSinCos)
$$
It follows easily that
$$
     \nu_k\le \nu_{k+1} \le \nu_k + 2, \qquad
     k=1,\dots,m-1.                                            \eqno(\eqExplOne)
$$

Let us denote
$$
    \beta_{k,0}^+ = A_k\cap\Big\{\,|\Arg\zeta-\theta_k|
    \le{\pi\ell_k^+\over\ell_k}\,\Big\},
    \qquad b_k^+ = \Area(\beta_{k,0}^+),
    \qquad k=1,\dots,m,
$$
where $\zeta = z/w$ is a standard complex coordinate on $\P^1$, and 
we shall choose the numbers
$\theta_1,\dots,\theta_m$ later.
In other words, the angular width of the domain $\beta_{k,0}^+$ is equal to 
$2\pi\ell^+_k/\ell_k = (2\pi s_k)/(n_ka_k)$. This implies
$$
     b_k^+ = {a_k\ell_k^+\over \ell_k} = {s_k\over n_k}.    \eqno(\eqExplTwo)
$$

Let us set
$$
     \beta_{k,j}^+ = R_{n_k}^j(\beta_{k,0}^+),  
                    \qquad
     \beta_{k,j}^- = \beta_{k,j}^+\cap\beta_{k,j+1}^+,
                    \qquad k=1,\dots,m, \quad j=1,\dots,n_k
$$
(here and further, when using the double index $(k,j)$, we assume that
$j$ is a residue mod $n_k$).

The angular width of $\beta_{k,j}^+$ is equal to $(2\pi s_k)/(n_ka_k)$ which
is not less than the angular size of the rotation $R_{n_k}$ (because $s_k/a_k\ge 1$).
Hence $\beta_{k,j}^- \ne\varnothing$ and
$$
     b_k^- := \Area(\beta_{k,j}^-) 
   = b_k^+ - {a_k\over n_k} = {s_k\over n_k} - {a_k\over n_k} = {s_{k-1}\over n_k}.
                                                            \eqno(\eqExplThree)
$$ 
Let $p_{k,j}$ be the midpoint of the arc $(\partial\Delta_k)\cap\beta_{k,j}^-$, and
let $q_{k,j}$ be the midpoint of the arc $(\partial\Delta_{k-1})\cap\beta_{k,j}^+$.
Now let us choose the numbers $\theta_k$ used in the definition of the domains 
$\beta_{k,0}^+$ so that $p_{k,0}=q_{k+1,0}$ for all $k$.
Since $p_{k,j} = R_{n_k}^j(p_{k,0})$ and $q_{k,j} = R_{n_k}^j(q_{k,0})$,
we have
$$
\xalignat2
      &\{p_{k,j}\,|\,0\le j<n_k\} = \{q_{k+1,j}\,|\,0\le j<n_{k+1}\}
        &&\text{ for $n_k=n_{k+1}$ }\\
      &\{p_{k,j}\,|\,0\le j<n_k\} \subset \{q_{k+1,j}\,|\,0\le j<n_{k+1}\}
        &&\text{ for $n_k<n_{k+1}$ }
\endxalignat
$$
Moreover, for all $k,j$ we have
$$
      p_{k,j} = q_{k+1,\mu_k j}, \qquad\text{where}\quad
      \mu_k = {n_{k+1}\over n_k} = 2^{\nu_{k+1} - \nu_k}.
$$

Note that by definition we also have
$$
     p_{m,1}=p_{m,2}=\dots=p_{m,n_m}=p_\infty,
     \qquad
     q_{1,1}=q_{1,2}=\dots=q_{1,n_1}=p_0.
$$
Let $\alpha_{k,j}^+$, $k=1,\dots,m$, $j=1,\dots,n_k$, be the path coming
from $p_{k,j}$ to $p_{k,j-1}$ along the boundary of
$\beta^+_{k,j}$ in the positive direction which passes any point of 
$(\partial\beta^+_{k,j})\setminus\{p_{k,j}\}$ at most once.
In the case $k=m$ this definition is ambiguous (because $p_{m,j}=p_{m,j-1}=p_\infty$),
but we assume that $\alpha_{m,j}^+$ is the complete loop around 
$\beta^+_{k,j}$ in the positive direction starting and finishing at $p_\infty$.
Let 
$\gamma_{k,j}^{(0)}$ (resp. $\gamma_{k,j}^{(1)}$) be the half of
$\alpha_{k,j}^+$ coming from $p_{k,j}$ to $q_{k,j}$
(resp. from $q_{k,j}$ to $p_{k,j-1}$).
Finally, let us set (see Figures {\figBetaPlus} and \figBetaMinus)
$$
      \alpha_{k,j}^- = \gamma_{k,j+1}^{(1)} + \gamma_{k,j}^{(0)},
      \qquad k=1,\dots,m, \quad j=1,\dots,n_k,
$$
$$
        \alpha_{m,1}=\alpha_{m,1}^+, \qquad
      \alpha_{k,1} =
           \alpha_{k,1}^+ - \sum_{j=0}^{\mu_k - 1}\alpha_{k+1,j}^-,
             \quad k=1,\dots,m-1,
$$
$$
      \alpha_{k,j+1} = R_{n_k}(\alpha_{k,j}),
      \qquad k=1,\dots,m, \quad j=1,\dots,n_k-1.
$$

\midinsert
\epsfxsize 125mm
\centerline{\epsfbox{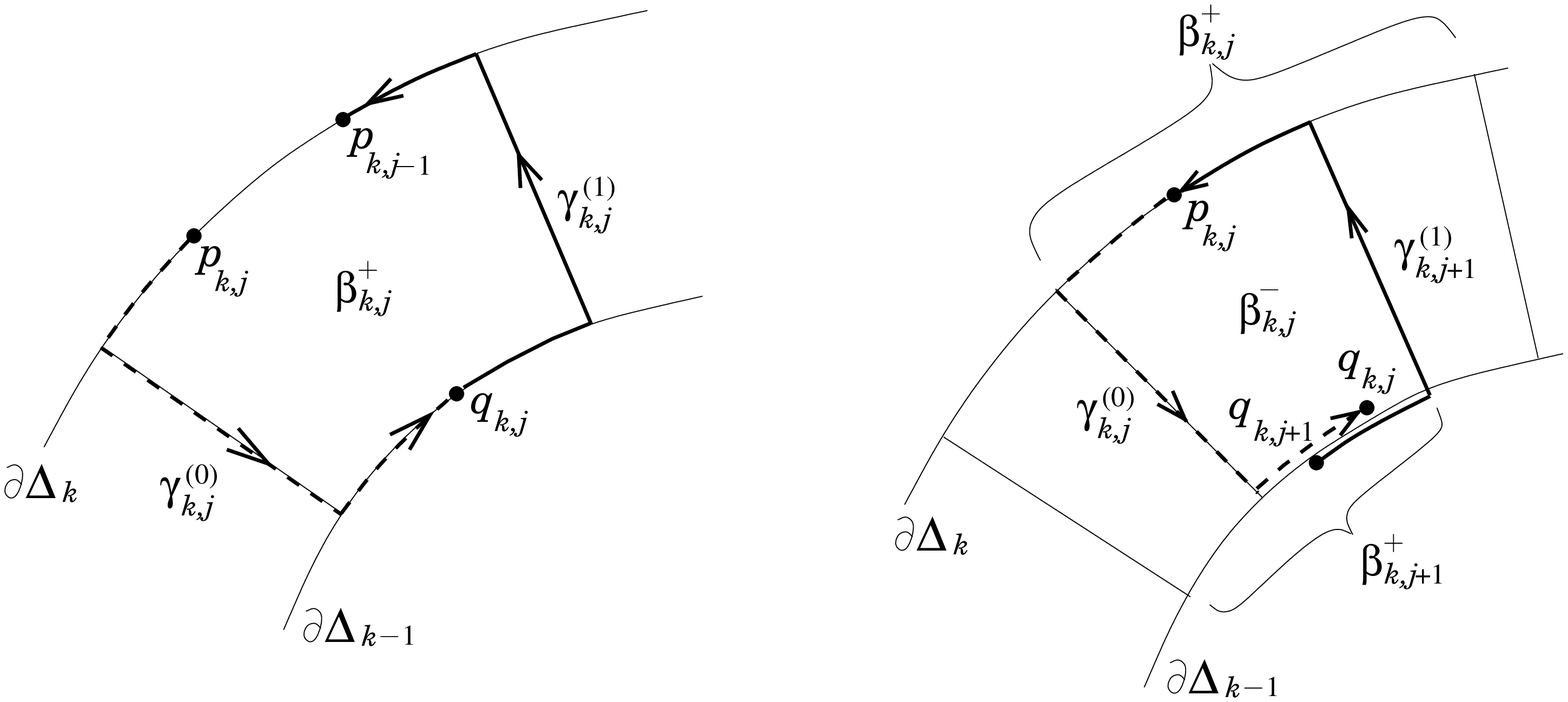}}
\botcaption{ Fig.~\figBetaPlus. 
     $\alpha_{k,j}^+ = \gamma_{k,j}^{(0)}+\gamma_{k,j}^{(1)}\,.$
        \hbox to 20mm{}
  Fig.~\figBetaMinus.
     $\alpha_{k,j}^- = \gamma_{k,j}^{(0)}+\gamma_{k,j+1}^{(1)}\,.$}
\endcaption
\endinsert

Let 
$\gamma:[0,1]\to\P^1$ be a piecewise smooth path and $\tilde p$ a point in $\SS^3$
such that $\pr(\tilde p)=\gamma(0)$ (as in \S\sectStokes, here $\pr$ denotes the
standard projection $\SS^3\to\P^1$). Then there exists a unique Legendrian path
$\tilde\gamma:[0,1]\to\SS^3$ such that
$\tilde\gamma(0)=\tilde p$ and $\pr\circ\tilde\gamma=\gamma$.
This follows from the fact that the fibers of $\pr:\SS^3\to\P^1$ are transverse
to the field of complex tangents  $\ker(\eta|_{\SS^3})$.
The path  $\tilde\gamma$ is called the {\it Legendrian lift} 
of $\gamma$ starting at $\tilde p$.

We shall construct Legendrian lifts 
$\tilde\alpha_{k,j}^\pm$ and $\tilde\alpha_{k,j}$
of 
$\alpha_{k,j}^\pm$ and $\alpha_{k,j}$ and we shall show that 
$\{\tilde\alpha_{k,j}\}$ is the required Legendrian net.

Let us set
$$
     \tilde p_{m,j} = ( e^{2\pi ij/n_m}, 0 )\in\C^2,
     \qquad j=1,\dots,n_m.
$$
The points $p_{m,j}$ belong to $\Gamma$, 
and we have $\tilde R_{n_m}(p_{m,j}) = p_{m,j+1}$.
Let $\tilde\gamma_{m,j}$, $j=1,\dots,n_m$, be the path
$[(j-1)/n_m,j/n_m]\to\Gamma$, $t\mapsto(e^{it},0)$. It goes along
$\Gamma$ from $\tilde p_{m,j-1}$ to $\tilde p_{m,j}$.
By (\eqExplTwo), we have $\Area(\beta_{m,j}^+)=b^+_m = s_m/n_m = \pi/n_m$.
Thus, 
$$
   \int_{\tilde\gamma_{m,j}}\eta = {1\over n_m}\int_\Gamma\eta = {2\pi\over n_m}
   = 2 b_m^+.
$$

The further construction will be recurrent (the induction by $k$).
Suppose that for some $k\le m$, we have constructed points
$\tilde p_{k,j}$ and paths $\tilde\gamma_{k,j}$ in $\SS^3$ such that
for any $j=1,\dots,n_k$ the following conditions hold
(as we have seen above, they do hold for $k=m$).
\roster
\item"(i)"
      $\partial\tilde\gamma_{k,j}=\tilde p_{k,j}-\tilde p_{k,j-1}$;
\item"(ii)"
      $\pr(\tilde p_{k,j}) = p_{k,j}$ and 
      $\pr(\tilde\gamma_{k,j})$ is the arc of $\partial\Delta_k$
      going in the positive direction from $p_{k,j-1}$ to $p_{k,j}$ and
      passing any point of $\partial\Delta_k$ at most once.
\item"(iii)"
      $\int_{\tilde\gamma_{k,j}}\eta = 2 b_k^+$;
\item"(iv)"
      $\tilde R_{n_k}( \tilde p_{k,j} ) = \tilde p_{k,j+1}$ and 
      $\tilde R_{n_k}( \tilde\gamma_{k,j} ) = \tilde\gamma_{k,j+1}$.
\endroster

Let $\tilde\alpha^+_{k,j}$ be the Legendrian lift of $\alpha^+_{k,j}$
starting at $\tilde p_{k,j}$.
Let us show that the end of $\tilde\alpha^+_{k,j}$ is $\tilde p_{k,j-1}$. 
Indeed, let us denote the end of $\tilde\alpha^+_{k,j}$ by $\tilde p$,
and let $[\tilde p,\tilde p_{k,j-1}]$ be the arc of the circle $\pr^{-1}(p_{k,j-1})$
from $\tilde p$ to $\tilde p_{k,j-1}$ chosen so that the cycle
$\tilde\alpha_{k,j}^+ + \tilde\gamma_{k,j} + [\tilde p,\tilde p_{k,j-1}]$
is zero-homologous in $\SS^3$. It is clear that the projection of this cycle to $\P^1$
coincides with $\partial\beta_{k,j}^+$.
Hence, Lemma {\lemStokes} implies
$$
   2b_k^+ = 2\Area(\beta^+_{k,j}) 
     = \int_{\tilde\alpha_{k,j}^+}\eta + \int_{\tilde\gamma_{k,j}}\eta
       + \int_{[\tilde p,\tilde p_{k,j-1}]}\eta
    = 0 + 2 b_k^+ + \int_{[\tilde p,\tilde p_{k,j-1}]}\eta.
$$
Therefore, $\int_{[\tilde p,\tilde p_{k,j-1}]}\eta=0$ and thus,
$\tilde p=\tilde p_{k,j-1}$.

Let us show that $\tilde R_{n_k}(\tilde\alpha_{k,j}^+)=\tilde\alpha^+_{k,j+1}$.
Indeed, let ${\Cal F_{k,j}}$ be the field of real tangent lines on the torus
$T_{k,j}=\pr^{-1}(\alpha^+_{k,j})$ which is cut by the field of complex tangents
$\ker\eta|_{\SS^3}$. Then 
$\tilde\alpha^+_{k,j}$ is the integral curve of $\Cal F_{k,j}$ passing through
$\tilde p_{k,j}$. It remains to note that $\tilde R_{n_k}$
takes $p_{k,j}$ and $T_{k,j}$ into $p_{k,j+1}$ and $T_{k,j+1}$ respectively.
Since, moreover, $\tilde R_{n_k}^*(\eta)=\eta$, it takes  $\Cal F_{k,j}$ into 
$\Cal F_{k,j+1}$.

Let $\tilde q_{k,j}$, $j=1,\dots,n_k$, be the point on $\tilde\alpha_{k,j}^+$
such that $\pr(\tilde q_{k,j}) = q_{k,j}$.
Let us set 
$\tilde\alpha_k = \sum_{j=1}^{n_k}\tilde\alpha^+_{k,j}$.
This is a closed spiral-like Legendrian curve on $\SS^3$ passing through the points
$\tilde p_{k,j}$ and $\tilde q_{k,j}$ and which is invariant under the rotation
$\tilde R_{n_k}$.
Let $\tilde\alpha_{k,j}^-$ be the Legendrian lift of $\alpha_{k,j}^-$ starting 
at $\tilde q_{k,j+1}$. Then
$\tilde\alpha_k = \sum_{j=1}^{n_k}\tilde\alpha^-_{k,j}$.
Moreover, the curve $\tilde\alpha_k$ is divided by the points
$\tilde p_{k,j}$ into the arcs $\tilde\alpha_{k,j}^+$, and it is divided 
by the points  $\tilde q_{k,j}$ into the arcs $\tilde\alpha_{k,j}^-$.

Let  $\gamma_{k,0}^-$ be the arc of $\partial\Delta_{k-1}$ going into the
positive direction from $q_{k,0}$ to $q_{k,1}$ and passing any point of
$\partial\Delta_{k-1}$ at most once. Let us choose a (non-Legendrian) lift
$\tilde\gamma_{k,0}^-$ of  $\gamma_{k,0}^-$ from
$\tilde q_{k,0}$ to $\tilde q_{k,1}$ so that the cycle
$\tilde\alpha^-_{k,0} + \tilde\gamma_{k,0}^-$ to be zero-homologous in
$\SS^3$. Its projection to $\P^1$ coincides with
$\partial\beta_{k,0}^-$, hence, by Lemma {\lemStokes} we have
$$
   2b_k^- = 2\Area(\beta^-_{k,0}) 
     = \int_{\tilde\alpha_{k,0}^-}\eta + \int_{\tilde\gamma_{k,0}^-}\eta
     = 0 + \int_{\tilde\gamma_{k,0}^-}\eta.
                                                       \eqno(\eqExplFour)
$$
Let us set $\tilde\gamma_{k,j}^- = \tilde R_{n_k}^j(\tilde\gamma_{k,0}^-)$,
$j=1,\dots,n_k$, and let
$$
   \tilde p_{k-1,0} = \tilde q_{k,0},
   \qquad
   \tilde\gamma_{k-1,1} = \sum_{j=0}^{\mu_{k-1}-1}\tilde\gamma_{k,j}^-,
$$
$$
   \tilde p_{k-1,j+1} = \tilde R_{n_{k-1}}(\tilde p_{k-1,j})
   \quad\text{ and }\quad
   \tilde\gamma_{k-1,j+1} = \tilde R_{n_{k-1}}(\tilde\gamma_{k-1,j}).
$$ 
To complete the recurrent construction, it remains to check that Conditions
(i)--(iv) are satisfied for the points $\tilde p_{k-1,j}$ and for the paths
$\tilde\gamma_{k-1,j}$.
Indeed, by (\eqExplFour) we have
$$
    \int_{\tilde\gamma_{k-1,1}}\eta 
   = \sum_{j=0}^{\mu_{k-1}-1}\int_{\tilde\gamma_{k,j}^-}\eta
   = 2\mu_{k-1} b_k^-.
$$
Combining this with (\eqExplThree) and $\mu_{k-1} = n_k/n_{k-1}$,
we obtain (iii) for $\tilde\gamma_{k-1,j}$.
The other conditions are evident.

\medskip

Finally, for $k=m\,$, let us set
$$
     \tilde\alpha_{m,j} = \tilde\alpha_{m,j}^{\leg} + \tilde\alpha_{m,j}^{\pt},
     \quad\text{ where }\quad
      \tilde\alpha_{m,j}^{\leg} = \tilde\alpha_{m,j}^+,\;\;
      \tilde\alpha_{m,j}^{\pt} = \tilde\gamma_{m,j},\quad
      j=1,\dots,n_m,
$$
and for $k<m\,$, let us set
$$
      \tilde\alpha_{k,1} = \tilde\alpha_{k,1}^{\leg} =
           \tilde\alpha_{k,1}^+ - \sum_{j=0}^{\mu_k - 1}\tilde\alpha_{k+1,j}^-,
$$
$$
      \tilde\alpha_{k,j+1} = \tilde\alpha_{k,j+1}^{\leg}
     = \tilde R_{n_k}(\tilde\alpha_{k,j}),
      \quad j=1,\dots,n_k-1.
$$
$$
      \tilde\alpha_{k,j}^{\pt}=0,\qquad
       j=1,\dots,n_k.
$$
(see Figure \figNet).

\midinsert
\epsfxsize 125mm
\centerline{\epsfbox{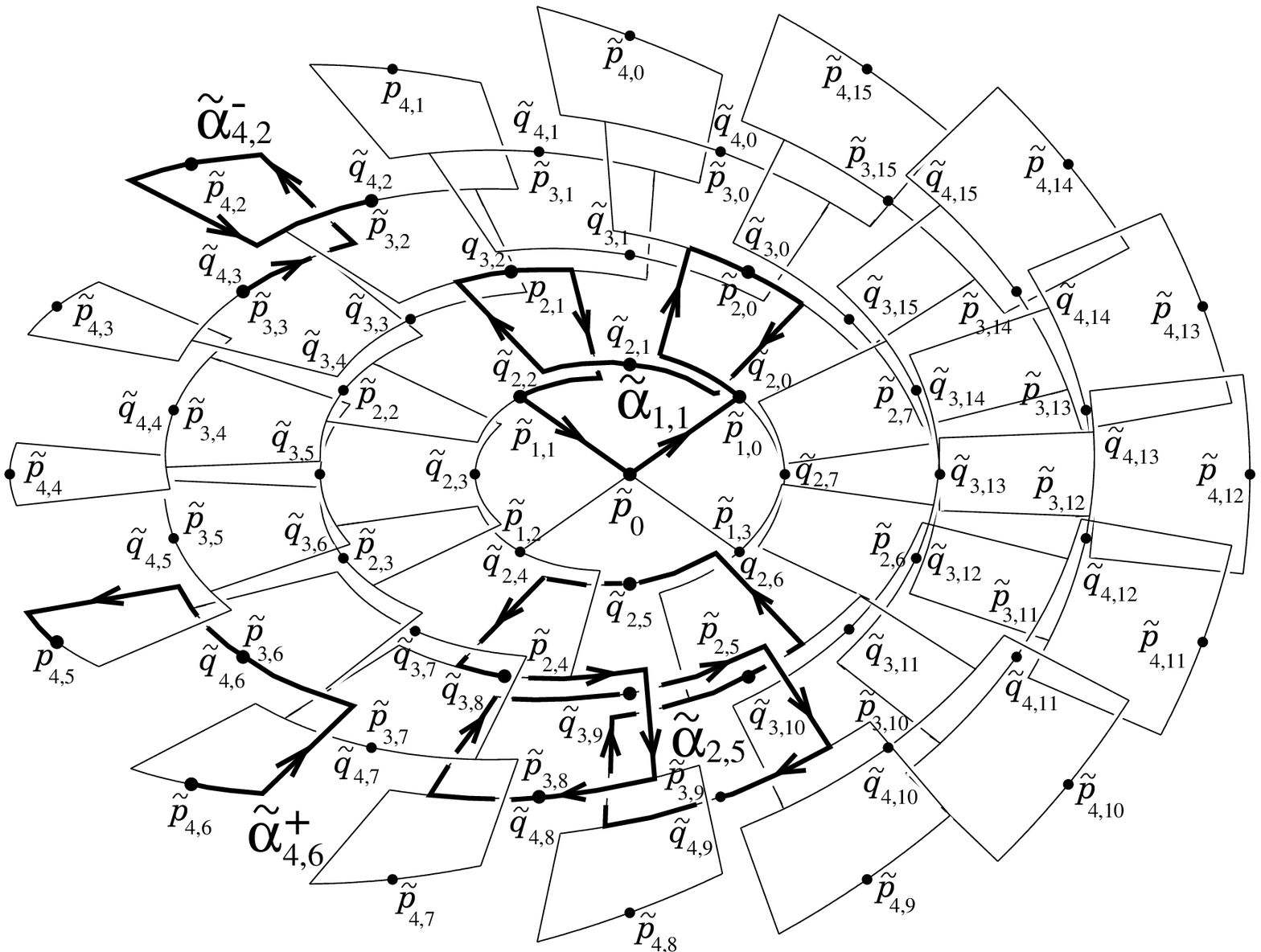}}
\botcaption{ Fig.~\figNet } 
     $\supp(\tilde\alpha_1+\dots+\tilde\alpha_4)$ for $(\nu_1,\dots,\nu_4)=(2,3,4,4)$.
\endcaption
\endinsert

Let us denote 
$\calA_\eps = \{\tilde\alpha_{k,j}\,|\,1\le k\le m,\, 1\le j\le n_k\,\}$.

\proclaim{ Proclaim \propExplOne }
        $\calA_\eps$ is a Legendrian net hanged on $\Gamma$.
\endproclaim

\demo{ Proof } It is easy to see that the chains
$\tilde\alpha_{k,j}$ are cycles.
Since the chains $\tilde\alpha_{k,j}^{\leg}$ 
(resp. $\tilde\alpha_{k,j}^{\pt}$) are Legendrian (resp. positively transverse)
by construction, it remains to check that  $\sum\tilde\alpha_{k,j}=\Gamma$.
Indeed, we have
$$
    \sum_{j=1}^{n_m}\tilde\alpha_{m,j} = 
    \sum_{j=1}^{n_m}\tilde\alpha_{m,j}^+ +\sum_{j=1}^{n_m}\tilde\gamma_{m,j}
   = \tilde\alpha_m + \Gamma,
$$
and for $k<m$, since $\tilde R_{n_k}=\tilde R_{n_{k+1}}^{\mu_k}$, we have
$$
   \tilde R_{n_k}^j(\tilde\alpha_{k+1,j'}^-)
    =\tilde R_{n_{k+1}}^{\mu_k j}(\tilde\alpha_{k+1,j'}^-)
    =\tilde\alpha_{k+1,j'+\mu_k j}^-, 
   \qquad\text{ hence, }
$$
$$
   \sum_{j=1}^{n_k}\tilde\alpha_{k,j} = 
    \sum_{j=1}^{n_k}\tilde\alpha_{k,j}^+ -
    \sum_{j=1}^{n_k}\sum_{j'=0}^{\mu_k - 1}\tilde R_{n_k}^j(\tilde\alpha_{k+1,j'}^-)
   = \sum_{j=1}^{n_k}\tilde\alpha_{k,j}^+ - \sum_{j=1}^{n_{k+1}}\tilde\alpha_{k+1,j}^-
   = \tilde\alpha_k - \tilde\alpha_{k+1}.
$$
Therefore,
$$
    \sum_{k,j}\tilde\alpha_{k,j} = (\tilde\alpha_1 - \tilde\alpha_2) +
    (\tilde\alpha_2 - \tilde\alpha_3) +\dots 
     + (\tilde\alpha_{m-1} - \tilde\alpha_m) + (\tilde\alpha_m + \Gamma) 
     = \tilde\alpha_1+\Gamma.
$$
It remains to note that $\tilde\alpha_1=0$ (see Figure \figNet) because
$\beta_{1,j}^-$ is a segment of a geodesic between 
$p_{0,j}=p_0$ and $p_{1,j}$, hence $\tilde\alpha_{1,j}^-=0$ for all $j=1,\dots,n_1$. 
\qed
\enddemo

\proclaim{ Proposition \propExplTwo}
$\length\tilde\alpha_{k,j}<\eps$ for all $k,j$.
\endproclaim

\demo{ Proof } It follows from
(\eqNorm) that the length of a path on $\P^1$ is equal to the length of its
Legendrian lift to  $\SS^3$. Hence,
$$
    \length\tilde\alpha_{k,j} = 
      \cases
             \length\partial\beta^+_{m,j} + \length\tilde\gamma_{m,j},
                     & k=m,\\
             \length\partial\beta^+_{k,j} + \mu_k\length\partial\beta^-_{k+1,j'}
                           & k<m.
      \endcases
$$
It is clear that
$$
   \length\partial\beta_{k,j}^\pm \le 2\ell_k^\pm 
   + 2\cdot(\text{width of $A_k$}) \le 2\ell_k^\pm + {\pi\over m} 
   \le 2\cdot{\eps\over 20} + {\eps\over 10} = {\eps\over 5},
$$
and (\eqExplOne) implies that $\mu_k\le 4$ for all $k$.
Hence, for $k<m$, we have
$$
   \length\tilde\alpha_{k,j} \le (1+\mu_k){\eps\over5} \le \eps. 
$$
It follows from (\eqEll) and (\eqSinCos) that
$$
  {1\over n_m}\le {\eps a_m\over20s_m\ell_m} =
    {\eps(1-\cos(\pi/m))\over 20\pi\sin(\pi/m)} = \eps\cdot O(1/m) = O(\eps^2).
$$
Thus,
$\length\tilde\alpha_{m,j} \le \length\partial\beta^+_{m,j} + O(\eps^2)
       \le (\eps/5) + O(\eps^2) \le\eps$.
\qed\enddemo

\proclaim{ Corollary \corExpl } 
An upper bound $n(\eps) \le \card\calA_\eps = O(\eps^{-3})$ holds.
\endproclaim

\demo{ Proof }
It follows from (\eqExplOne) that $n_1\le n_2\le\dots\le n_m$, 
hence
$\card\calA_\eps \le m n_m$.
It is clear that $m=O(\eps)$, and it easily follows from (\eqSinCos) 
that $n_m=O(\eps^{-2})$.
Therefore, $\card\calA_\eps = O(\eps^{-3})$.

The estimate $n(\eps) \le \card\calA_\eps$ follows from the construction given in
\S\sectAppr. 
\qed\enddemo

\remark{ Remark \remExpl }
By Proposition \propGeneric, it is not important for us if the set 
$\calA_\eps$ is generic or not. However, it is such everywhere except
the point $p_0$ (see Figure \figNet). It one changes slightly the parameters 
of the construction of $\calA_\eps$, it is not difficult to achieve
$n_1\le3$. In this case, $\calA_\eps$ will be generic everywhere, including $p_0$.
\endremark


\enddocument